\documentclass[
12pt,
DIV=20, 
twoside, 
headings=small, 
BCOR1cm, 
toc=graduated,
]
{scrartcl}

\usepackage{enumerate}
\usepackage{amsmath} 
\usepackage{amssymb}	
\usepackage{hyperref} 

\usepackage[numbers,sort&compress]{natbib}

\newcommand{\arpion}[1]{\footnote{#1}}

\newcommand{\Cov}{{\rm Cov}}
\newcommand{\D}{\mathbb{D}}
\newcommand{\dconv}{\overset{\mathcal{D}}{\longrightarrow}}
\newcommand{\dequal}{\overset{\mathcal{D}}{=}}
\newcommand{\E}{{\rm E}}
\newcommand{\F}{F}
\newcommand{\goodcl}{\mathcal{G}}
\newcommand{\gu}{\ra\infty} 

\newcommand{\lgauss}[1]{\left\lfloor #1 \right\rfloor}
\newcommand{\linv}{^{\leftarrow}}
\newcommand{\lra}{\longrightarrow}

\newcommand{\N}{\mathbb{N^*}}
\newcommand{\ngu}{n\ra\infty} 
\newcommand{\No}{\mathbb{N}}

\newcommand{\R}{\mathbb{R}}
\newcommand{\ra}{\rightarrow}

\newcommand{\rinv}{^{\rightarrow}}

\newcommand{\ugauss}[1]{\left\lceil #1 \right\rceil}

\newcommand{\mult}[1]{\overline{#1}}

\newcommand{\Z}{\mathbb{Z}}
\newcommand{\1}{\mathbf{1}}

\newcommand{\balance}[2]{{$#1$}-{$#2$}-control}
\newcommand{\MMP}[3]{\mathbf{MM}_{#1,#2}(#3)}

\renewcommand{\epsilon}{\varepsilon}

\renewcommand{\L}{L}
\renewcommand{\P}{{\rm P}}
\renewcommand{\phi}{\varphi}
\renewcommand{\theta}{\vartheta}

\usepackage{amsthm} 

	\newtheorem{theorem}{Theorem}
	\newtheorem{proposition}{Proposition}
	\newtheorem{lemma}{Lemma}
	\newtheorem{corollary}{Corollary}

\theoremstyle{definition}

	\newtheorem{definition}{Definition}
	\newtheorem{example}{Example}
	\newtheorem{remark}{Remark}

	\usepackage{color}
\newcommand{\CUT}[1]{}
\newcommand{\CUTT}[1]{}

\newcommand{\change}{}

\begin{document}


\title{An Empirical Process Central Limit Theorem for Multidimensional Dependent Data%
\footnote{
The final publication is available at \newline {http://www.springerlink.com/openurl.asp?genre=article\&id=doi:10.1007/s10959-012-0450-3}
}
}

\author{Olivier Durieu\footnote{Laboratoire de Math\'ematiques et Physique Th\'eorique (UMR CNRS 7350), F\'ed\'eration Denis Poisson (FR CNRS 2964), Universit\'e de Tours, Parc de Grandmont, 37200 Tours, France,
olivier.durieu@lmpt.univ-tours.fr, Fax: (+33) 247 367 068, Fon:  (+33) 247 367 421 
}
, Marco Tusche\footnote{Fakult\"at f\"ur Mathematik, Ruhr-Universit\"at Bochum, Universit\"atsstra\ss e 150, 44780 Bochum, Germany,
marco.tusche@RUB.de, Fax: (+49) 234 32 14039, Fon: (+49)234 32 23425
}
}
\date{October 1., 2012}
\maketitle



\begin{abstract}\noindent{\sffamily\bfseries Abstract}
Let $(U_n(t))_{t\in\R^d}$ be the empirical process associated to an $\R^d$-valued stationary process $(X_i)_{i\ge 0}$.
In the present paper, we introduce very general conditions for weak convergence of $(U_n(t))_{t\in\R^d}$, which only involve properties of processes $(f(X_i))_{i\ge 0}$ for a restricted class of functions $f\in\goodcl$. 
Our results significantly improve those of \cite{DehDurVol09} and \cite{DehDur11} and provide new applications.

The central interest in our approach is, that it does not need the indicator functions, which define the empirical process $(U_n(t))_{t\in\R^d}$, itself to belong to the class $\goodcl$.
This is particularly useful when dealing with data arising from dynamical systems or functionals of Markov chains.
In the proofs we make use of a new application of a chaining argument and generalize ideas first introduced in \cite{DehDurVol09} and \cite{DehDur11}.

Finally we will show how our general conditions apply in the case of multiple mixing processes of polynomial decrease and causal functions of independent and identically distributed
processes, which could not be treated by the preceding results in \cite{DehDurVol09} and \cite{DehDur11}.
\end{abstract}


\section{Introduction}

The present paper concerns the question of the weak convergence of empirical processes under weak dependence of the underlying process.
Let us consider some stationary $\R^d$-valued process\arpion{$\No:=\{0,1,\ldots\},$ $\N:=\No\backslash\{0\}$} $(X_i)_{i\in\No}$. 
The empirical process of $(X_i)_{i\in\No}$ is the process $(U_n)_{n\in\N}$ of $\D([-\infty,\infty]^d)$-valued\arpion{\label{txt:multi_cadlag}Let $\D([-\infty,\infty]^d)$ be the space of generalized multidimensional c\`adl\`ag functions $[-\infty,\infty]^d\ra\R$ (for definition see \cite[p.1286]{Neu71}), equipped with the multidimensional Skorohod metric $d_0$ as introduced in \cite[p.1289]{Neu71} (see also \cite{Str72}). 
Note that $(\D([-\infty,\infty]^d),d_0)$ is a complete and separable space 
(more precisely, \cite{Neu71} and \cite{Str72} proved this for the space $\D([0,1]^d)$, but -- since $[0,1]$ and $[-\infty,\infty]$ are homeomorphic -- the metric on $\D([-\infty,\infty]^d)$ can naturally be extended to a metric on $\D([-\infty,\infty]^d)$ which conserves all relevant properties (c.f. \cite[p.1081f.]{DehDur11})).} random variables 
$$U_n:=(U_n(t_1,\ldots,t_d))_{(t_1,\ldots,t_d)\in [-\infty,\infty]^d}$$
given by
\begin{align*} 
	U_n(t_1,\ldots,t_d) :=& \frac{1}{\sqrt{n}}\sum_{i=1}^n \Bigl( \1_{\prod_{j=1}^d[-\infty,t_j]}(X_i) - \E\1_{\prod_{j=1}^d[-\infty,t_j]}(X_0) \Bigr). 
\end{align*}


The study of the weak convergence of such a process began with Donsker \cite{Don52}, who proved convergence of the empirical process to some Gaussian process in the case of i.i.d.\ $\R$-valued data.
Donsker's result has been extended to sequences of weakly dependent ($\R$ or $\R^d$-valued) random variables by many authors.
Among others, it shall be remarked that 
Billingsley \cite{Bil68} gave a result for functionals of $\phi$-mixing sequences, Berkes and Philipp \cite{BerPhi77} under strong mixing assumptions, Doukhan, Massart and Rio \cite{DouMasRio95} for absolutely regular sequences, Borovkova, Burton and Dehling \cite{BorBurDeh01} for functionals of absolutely regular processes,
Doukhan and Louhichi \cite{DouLou99}, and Dedecker and Prieur \cite{DedPri07} for new dependence coefficients.

The technique of the proof usually consists of establishing finite-dimensional convergence of the process
and to use an appropriate moment bound to get tightness of the process (see \cite{DehPhi02}).

In the recent paper \cite{DehDurVol09}, the authors introduced a new approach which is useful when the required properties on the underlying process are only established for a class of functions not containing the indicators. The required conditions are the Central Limit Theorem and a bound of the $4$-th moment of partial sums for the process $(f(X_i))_{i\ge 0}$, when $f$ belongs to the class of Lipschitz functions. Under these conditions, in the case of $\R$-valued data, they were able to prove an empirical CLT. Later, in \cite{DehDur11}, the technique was adapted to treat the case of $\R^d$-valued data and using the class of H\"older functions.
In that case bounds on higher moments are needed.

These required conditions can be satisfied in cases where the data arise from a dynamical system. In particular, for a large class of dynamical systems, they can be derived from the study of the so called Perron-Frobenius operator. More generally, they can be established in cases where the system presents some multiple mixing properties (see Section \ref{sec:main-result} and \cite{DehDur11} for more details).

The aim of the present paper is to give a very general setting under which the technique of \cite{DehDurVol09} and \cite{DehDur11} can be applied.
In this way, we obtain extensions of the previous results. In particular, we show that the approach of \cite{DehDurVol09} and \cite{DehDur11} can be applied with more general classes of functions than Lipschitz and H\"older functions. We provide general conditions on the class of functions for which the CLT and moment bounds have to hold. In addition, we use a more general moment bound, balancing conditions on the distribution of $X_0$ with properties of the class of functions.

Our results allow us to handle processes with weaker moment bound than in previous works. In particular, we are able to treat systems satisfying a multiple mixing property with polynomial rates, whereas in \cite{DehDur11} exponential rates were required. As a concrete application, we obtain new result for causal functions of i.i.d.\ processes.

\section{Central Conditions and Statement of Main Results}\label{sec:main-result}

\subsection{Conditions}\label{sec:conditions}
Before coming to the statement of our main results, let us have a look at our conditions, what they are needed for and some example of how they can be established.

As pointed out in the introduction, the present paper presents a technique which differs from the usual finite-dimensional convergence plus tightness of $U_n$ approach.
This is useful in cases, where one cannot or not without great effort directly prove these conditions, but where one can establish similar results for functions $(f(X_i))_{i\in\No}$ of the underlying process. 
These situations appear, e.g.\ when dealing with data arising from Markov chains or dynamical systems (c.f.\ \cite{DehDurVol09} and \cite{DehDur11}).

\subsubsection{Central Limit Theorem under a Class of Functions.}

We say that a process $(X_i)_{i\in\No}$ satisfies a Central Limit Theorem under a class of functions $\mathcal{G}$, if for every $f\in\mathcal{G}$ such that $\E(f(X_0))=0$, there exists $\sigma^2_f<\infty$ such that
\begin{align} 
		\frac{1}{\sqrt{n}}\sum_{i=1}^n f(X_i) \dconv N(0,\sigma^2_f).  \label{con:clt}
	\end{align}

A lot of research has been devoted to establish CLTs under various classes of functions, such as functions of bounded variation, Lipschitz functions and H\"older functions. Hennion and Herv\'e \cite{HenHer01} give a survey of such results for dynamical systems and Markov chains, where the Perron-Frobenius operator or the Markov operator satisfies a spectral gap condition. Another example are ergodic torus automorphisms, for which Leonov \cite{Leo60} and Le Borgne \cite{Leb99} established CLTs for H\"older functions. In the present paper, we consider causal functions of i.i.d.\ processes where the CLT holds for H\"older functions. CLTs for Lipschitz functions were studied for random iterative Lipschitz models and linear function of i.i.d.\ processes in \cite{DehDurVol09}.

\subsubsection{Moment Bounds under a Class of Functions}\label{sec:moment-bounds}

To prove tightness, one usually makes use of conditions on the moments of increments of the empirical process. 
In order to restrict our conditions to functions of the process we will work with the following type of moment bounds under the same function space $\mathcal{G}$ for which the CLT \eqref{con:clt} holds:

There are finite constants $C>0$, $r\ge 1$, $p\in\N$ and nondecreasing functions $\Phi_1,\ldots,\Phi_p : \R_0^+\ra\R_0^+$ such that for all $f\in\mathcal{G}$
with $\|f\|_{\infty} \leq 1$ and all $n\in\N$ we have
	\begin{align}	
		&\E\biggl( \Bigl| \sum_{i=1}^n f(X_i) -\E(f(X_0)) \Bigr|^{2p} \biggr) \notag\\
		&\leq C \sum_{i=1}^p n^i \: \| f(X_0) - \E(f(X_0)) \|_r^i \: {\Phi_i} (\|f -\E(f(X_0))\|_{\goodcl}) 
		\label{con:2p-moment}.
\end{align}	
This condition is met, e.g. for processes satisfying a multiple mixing property.


\begin{definition}[Multiple mixing property]
Let $(X_i)_{i\in\No}$ be a stationary stochastic process of $\R^d$-valued random variables and let $\mathcal{G}$ be a class of measurable real-valued functions defined on $\R^d$ and equipped with a seminorm $\|\cdot\|_{\mathcal{G}}$.
For integers $i_1,...,i_j$, we write $i_j^*:=i_1+ \ldots + i_j$.
We say that $(X_n)_{n\in\No}$ 
\change is $\Theta,r$-multiple mixing 
with respect to $\mathcal{G}$
if there exist a constant $r\in[1,\infty)$ and a nonincreasing function 
$\Theta:\No\lra\R_0^+$ such that for any $p\in\N$ there is a constant $K_p<\infty$ satisfying:
\begin{align}
&\bigl| \Cov\bigl(f(X_0)f(X_{i_1^{*}}) \cdot\ldots\cdot f(X_{i_{q-1}^{*}})\,,\,f(X_{i_{q}^{*}})f(X_{i_{q+1}^{*}})\cdot\ldots\cdot f(X_{i_{p}^{*}}) \bigr) \bigr| \notag\\
&\leq K_p\|f(X_0)\|_r \|f\|_{\mathcal{G}} \Theta(i_q)
\label{con:multmix_pol}
\end{align}
for all $f\in\mathcal{G}$ with $\|f\|_{\infty}\leq 1$ and $\E(f(X_0))=0$ all $i_1,\ldots,i_p\in\No$, $q\in\{1,\ldots,p\}$.
In that case we write $(X_n)_{n\in\No}\in\MMP{\Theta}{r}{\mathcal{G}}$. 
\end{definition} 
The following proposition shows that multiple mixing systems satisfy a $2p$-th moment bound \eqref{con:2p-moment}.

\begin{proposition}\label{pro:multimix->2p-moment}
Let $(X_n)_{n\in\No}\in\MMP{\Theta}{r}{\mathcal{G}}$ for a function $\Theta$ such that there exists $p\in\N$ satisfying 
\begin{align}
	\sum_{i=0}^{\infty} i^{2p-2}\Theta(i) < \infty. \label{eq:temp020}
\end{align}
Then there is a $C>0$ such that \eqref{con:2p-moment} holds
for all $f\in\mathcal{G}$ such that $\|f\|_\infty\le 1$, with $r$ and $p$ as
above and $\Phi_i(x)=x^i$.
\end{proposition}
The proof of Proposition \ref{pro:multimix->2p-moment} will be given in Section \ref{sec:new_technique_applications_proof}.

Note that similar moment bounds have been obtained in \cite{DehDur11} under the
stronger assumption of multiple mixing with exponential rate, i.e. when
$\Theta(i)$ decreases exponentially. Under this assumption, one can get sharper
moment bounds, with $\Phi_i(x)=\log^{2p-i}(x+1)$. The moment bounds obtained in
the present paper are sufficient to apply Theorem \ref{the:newtech}, and thus
our generalization allows one to treat processes that are multiple mixing with
polynomial rate.


\subsubsection{Approximation of the Indicator Functions}
Conditions \eqref{con:clt} and \eqref{con:2p-moment} refer to the processes $(f(X_i))_{i\in\No}$ for $f\in\mathcal{G}$. In order to obtain results for the empirical process $U_n$, we need to approximate the indicator functions occurring in the definition of $U_n$ by functions from $\mathcal{G}$. To describe the quality of this approximation with respect to $\|\cdot\|_\mathcal{G}$, we introduce the upcoming definition.

In the following, $\leq$, $<$, ...\ used in $\R^d$ are to be understood component-wise\arpion{i.e.\ for $a=(a_1,\ldots,a_d)$, $b=(b_1,\ldots,b_d)\in [-\infty,\infty]^d$ write $ a \leq b$ if and only if $a_i\leq b_i$ for all $i=1,\ldots,d$} 
\change
and sets such as $\{x\in [-\infty,\infty]^d : a \leq x < b \}$ are denoted by $(a,b]$.
We will also use the notation $w_g$ for the modulus of continuity of a real-valued function $g$, which is defined by
\begin{align*}
	w_{g}(\delta)&:=\sup \{ |g(t)-g(s)|: s,t\in\R^d,\ \|t-s\|\leq \delta\}.
\end{align*}

\begin{definition}[\balance{\mathcal{G}}{\F} function]\label{def:control}
Let $\mathcal{G}$ be some vector space of real-valued functions defined on $\R^d$, equipped with a seminorm $\|\cdot\|_{\mathcal{G}}$, and let $\F$ denote the multidimensional distribution function of an $\R^d$-valued random variable $X$.
We call a nondecreasing function $\Psi:\R_0^+\ra\R_0^+$ a \balance{\mathcal{G}}{\F}, if 
for every $a < b\in[-\infty,\infty]^d$ there is a function $\phi_{(a,b)}\in\mathcal{G}$ such that for any $x\in\R^d$
\change	\begin{align}
		\1_{(-\infty,a]} \leq \phi_{(a,b]} \leq \1_{(-\infty,b]}. \label{con:existence_phi}
	\end{align}
and such that 
		\begin{align}
			\| \phi_{(a,b)} \|_{\mathcal{G}}  \leq \Psi\biggl( \frac{1}{\underset{i=1,\ldots,d}{\min} w_{\F_i} (b_i-a_i)} \biggr) \label{con:normbounding},
		\end{align}
where $\F_i$ denotes the $i$-th marginal distribution function of $X$.

If such a function $\Psi$ exists, we say that $\mathcal{G}$ approximates the indicator functions (of rectangles $[-\infty,t]$, $t\in\R^d$) with \balance{\mathcal{G}}{\F} $\Psi$.
\end{definition} 

\begin{example}
As an example consider the class of bounded $\alpha$-H\"older functions $\mathcal{H}_{\alpha}$
equipped with the $\alpha$-H\"older norm
\[ \|f\|_{\mathcal{H}_{\alpha}}:= \|f\|_\infty + \sup \biggl\{ \frac{|f(x)-f(y)|}{\|x-y\|^{\alpha}} : {x,y\in\R^d,\ x\neq y} \biggr\}. \]
Choose $\phi_{(a,b)}\in\mathcal{H}_{\alpha}$ as
\begin{align*}
	\phi_{(a,b)}(x_1,\dots,x_d):=\prod_{i=1}^d \phi \Bigl( \1_{(-\infty,\infty)^2}(a_i,b_i) \cdot \frac{x_i - b_i}{b_i-a_i} \Bigr),
\end{align*}
where $\phi:\R\ra\R$ is given by $\phi(x)=\1_{[-\infty,-1]}(x)-x\1_{(-1,0]}(x)$.
Obviously this choice of $\phi_{(a,b)}$ satisfies \eqref{con:existence_phi}.
Let us now check condition \eqref{con:normbounding}.
Since for all $j=1,\ldots,d$, $w_{\F_j}(\delta)\le w_\F(\delta)$, we have
	\begin{align*}
		b_j-a_j 
		&\geq \inf\bigl\{ \delta>0 : w_{\F}(\delta) \geq \min_{i=1,\ldots,d} w_{\F_i}(b_i-a_i) \bigr\}.
	\end{align*}
Thus, by the definition of $\phi_{(a,b)}$ we obtain
\begin{align*}
	\|\phi_{(a,b)}\|_{\mathcal{H}_{\alpha}} 
	& \leq d \max_{j=1,\ldots,d} \1_{(-\infty,\infty)^2}(a_j,b_j) \cdot \frac{1}{(b_j-a_j)^{\alpha}} +1 \\
	& \leq d \Bigl( w_{\F}\linv \Bigl({\min_{i=1,\ldots,d}  w_{\F_i}(b_i-a_i)}\Bigr)\Bigr)^{-\alpha} +1,
\end{align*}
where $w_{\F}\linv(y):=\inf \bigl\{\delta>0 : w_{\F}(\delta) \geq {y}\bigr\}$.
Hence \eqref{con:normbounding} is satisfied for the nondecreasing function $\Psi$ given by
\begin{align} 
	\Psi(z) := d \Bigl( w_{\F}\linv \Bigl(\frac{1}{z}\Bigr)\Bigr)^{-\alpha} +1 
	\label{eq:control}
\end{align}
and thus $\Psi$ defines an \balance{\mathcal{H}_{\alpha}}{\F}, which gives us the following lemma:
\begin{lemma}\label{lem:control}
The space of bounded $\alpha$-H\"older functions $\mathcal{H}_{\alpha}$ approximate the indicator functions with \balance{\mathcal{H}_{\alpha}}{\F} $\Psi$ given by \eqref{eq:control}.
\end{lemma}

\end{example}

\subsection{Main Theorems}

Our main result is the following:

\begin{theorem}\label{the:newtech} Let $(X_i)_{i\in\No}$ be a stationary process of $\R^d$-valued random vectors with continuous multidimensional distribution function $\F$.
\noindent Assume that there is a vector space $\mathcal{G}$ of measurable functions $\R^d\ra\R$, containing the constant functions, equipped with a seminorm $\|\cdot\|_{\mathcal{G}}$, satisfying the following conditions:
\begin{enumerate}[(i)]
	\item For every $f\in\mathcal{G}$ such that $\|f\|_{\infty}<\infty$, the CLT \eqref{con:clt} holds.
	\item $\mathcal{G}$ approximates the indicator functions with \balance{\mathcal{G}}{\F} $\Psi$.
	\item There are constants $r\geq 1$,  $p>rd$, $\gamma_1,\ldots,\gamma_p$ satisfying
\begin{align}
	0\le \gamma_i< \frac{i}{r} + 2(p-i) - d\ \ \text{ for all } i=1,\ldots,p \label{con:gamma_i},
\end{align}
and some nondecreasing functions $\Phi_1,\ldots,\Phi_p: \R_0^+\ra\R_0^+$ satisfying
	\begin{align} 
		{\Phi_i} (2\Psi(z)) = \mathcal{O} (z^{\gamma_i})\ \ \text{ as }z\gu
		\label{con:Phi_of_Psi},
	\end{align}	
such that for every
$f\in\mathcal{G}$ with $\|f\|_{\infty}\leq1$ the moment bound \eqref{con:2p-moment} holds. 
\end{enumerate}
Then there is a centered Gaussian process $(W(t))_{t\in[-\infty,\infty]^d}$ with almost surely continuous sample paths such that $U_n \dconv W$, in the space $\D([-\infty,\infty]^d)$.
\end{theorem}

\begin{remark}
 In fact, the functions $\Phi_1,\ldots,\Phi_p$ and the function $\Psi$ only have to be nondecreasing for sufficiently large arguments.
 Note that condition \eqref{con:2p-moment} has only to be satisfied for a certain subclass of $\mathcal{G}$, see Remark \ref{rem:class_G}.
\end{remark}

As a \change consequence of this abstract theorem, we can give a statement for multiple mixing processes, for which conditions are more easily verifiable.

\begin{theorem}\label{the:multimix}
Let $(X_i)_{i\in\No}$ be a stationary $\R^d$-valued process with continuous multidimensional distribution function $\F$. Assume there is a vector space $\mathcal{G}$ of measurable functions $\R^d\ra\R$, containing the constant functions, equipped with a seminorm $\|\cdot\|_{\mathcal{G}}$that  satisfies the following conditions:
\begin{enumerate}[(i)]
\item\label{2-i} For every $f\in\mathcal{G}$ such that $\|f\|_{\infty}<\infty$ the CLT \eqref{con:clt} holds.
	\item\change\label{2-ii} 
The process $(X_i)_{i\in\N}$ is $\Theta,r$-multiple mixing with respect to $\mathcal{G}$ for some $r\geq 1$ and a $\Theta:\No \ra \R^+_0$, such that there exists a $p>dr$ satisfying $\sum_{i=0}^{\infty} i^{2p-2}\Theta(i)<\infty$.
\item $\mathcal{G}$ approximates the indicator functions with \balance{\mathcal{G}}{\F} $\Psi$ such that $\Psi(z)=\mathcal{O}(z^{\frac{1}{\gamma}})$ for some $\gamma > \frac{rp}{p-rd}$.
\end{enumerate}
Then there is a centered Gaussian process $(W(t))_{t\in[-\infty,\infty]^d}$ with almost surely continuous sample paths such that $U_n \dconv W$, in the space $\D([-\infty,\infty]^d)$.
\end{theorem}

\proof
By Proposition \ref{pro:multimix->2p-moment}, \eqref{con:2p-moment} holds with $\Phi_i(x)=x^i$. Then, taking $\gamma_i=\frac{i}{\gamma}$, the assumptions of Theorem \ref{the:newtech} are satisfied.
\qed

The rest of the paper is organized as follows: In Section \ref{sec:new_technique_applications} we present some particular cases of the theorem and some applications. 
The proof of Theorem \ref{the:newtech} is given in Section \ref{sec:new_technique_proof}.
Section \ref{sec:new_technique_applications_proof} is devoted to the proofs of the results of Section \ref{sec:new_technique_applications}.

\section{Applications}\label{sec:new_technique_applications}

While Theorem \ref{the:newtech} is stated in a very general context, this section is devoted to the special case where $\mathcal{G}=\mathcal{H}_{\alpha}$, the class of bounded $\alpha$-H\"older functions.
We will first obtain a corollary to Theorem \ref{the:newtech} specifying the conditions for this choice of $\mathcal{G}$. As an application, we will study empirical processes of causal functions of an i.i.d.\ process.

\subsection{Assumptions on H\"older continuous functions}

The case of H\"older continuous functions has already been considered in \cite{DehDur11}. In the present paper, we obtain generalizations of the results in \cite{DehDur11}, by allowing a larger class of functions $\Phi_i$ in the moment bound \eqref{con:2p-moment}. In this way, we are able to treat processes that satisfy a multiple mixing condition with polynomial decay, rather than the exponential decay considered in \cite{DehDur11}.

\begin{theorem}\label{pro:newtech_hoelder}
Let $(X_i)_{i\in\No}$ be a stationary $\R^d$-valued process with continuous multidimensional distribution function $\F$.
Assume that:
\begin{enumerate}[(i)] 
	\item For every $f\in\mathcal{H}_{\alpha}$ the CLT \eqref{con:clt} holds.
	\item\label{3-ii} There are constants $\beta>0$, $\gamma>1$, $r\geq 1$, $p>rd$, $\kappa_1,\ldots,\kappa_p>0$, $\lambda>1$, $z_0\in\R^+$ and
an invertible increasing function $\Phi:\R_0^+\lra\R_0^+$ satisfying
\begin{align} 
&	w_{\F} (y) \leq \beta \Phi(y^{-{\alpha}})^{-{\gamma}} \label{con:momentbound'},
\end{align}
\begin{align}
 &\kappa_i < \gamma \Bigl( \frac{i}{r} + 2(p-i) -d \Bigr),\ \ i=1,\ldots,p, \label{con:kappa_i}
\end{align}
\begin{align} 
&	\Phi(2 z) \leq \lambda \Phi(z) \label{eq:phi_hoelder}, \ \ \forall z\geq z_0,
 \end{align}
such that for all $f\in\mathcal{H}_{\alpha}$ with $\|f\|_{\infty}\leq1$, the moment bound \eqref{con:2p-moment} holds with $\Phi_i=\Phi^{\kappa_i}$.
	\end{enumerate}
Then there is a centered Gaussian process $(W(t))_{t\in[-\infty,\infty]^d}$ with almost surely continuous sample paths such that $U_n \dconv W$, in the space $\D([-\infty,\infty]^d)$.
\end{theorem}

\proof
Using Lemma \ref{lem:control}, it is sufficient to show condition \eqref{con:Phi_of_Psi} of Theorem \ref{the:newtech} for the function
$\Psi$ defined in \eqref{eq:control} and for suitable $\gamma_1,\ldots,\gamma_p$ satisfying \eqref{con:gamma_i}.

By \eqref{con:momentbound'}, for $z>0$ we have
\[ 
w_{\F}\linv \Bigl(\frac{1}{z}\Bigr) \geq \inf \Bigl\{ x > 0 : \beta (\Phi(x^{-{\alpha}}))^{-{\gamma}} \geq \frac{1}{z} \Bigr\} 
 = \inf\biggl\{ x > 0 : \Phi(x^{-{\alpha}}) \leq \Bigl({\beta z}\Bigr)^{\frac{1}{\gamma}} \biggr\} 
\]
and therefore, since $\Phi$ is invertible, we obtain
\[ 
	w_{\F}\linv \Bigl(\frac{1}{z}\Bigr)^{-\alpha} 
\geq \inf\Bigl\{ x > 0 : x \geq \bigl(\Phi^{-1} \bigl(({\beta z}
)^{\frac{1}{\gamma}}\bigr)\bigr)^{-\frac{1}{\alpha}} \Bigr\}^{-\alpha} 
= \Phi^{-1} \Bigl(({\beta z})^{\frac{1}{\gamma}}\Bigr).
\]
Using \eqref{eq:phi_hoelder}, this yields that there is a finite constant $C'\in\R$ such that
\begin{align*} 
	\Phi_i(2\Psi(z))&=\Phi_i\Bigl(2d \Bigl( w_{\F}\linv \Bigl(\frac{1}{z}\Bigr)\Bigr)^{-\alpha} +2\Bigr) \notag\\
	&\leq \Bigl(\Phi(2d\Phi^{-1}((\beta z)^{\frac{1}{\gamma}}) +2)\Bigr)^{\kappa_i} \notag\\
	&\leq C' \Bigl(\Phi(\Phi^{-1}((\beta z)^{\frac{1}{\gamma}}))\Bigr)^{\kappa_i} \notag \\
	&= C'\beta^{\frac{\kappa_i}{\gamma}} \cdot z^{\frac{\kappa_i}{\gamma}}  
\end{align*}
for every $i\in\{1,\ldots,p\}$ and all sufficiently large $z$. Here, by \eqref{con:kappa_i}, each exponent $\gamma_i=\kappa_i/\gamma$ is positive and strictly smaller than $i/r+2(p-i)-d$. Thus conditions \eqref{con:Phi_of_Psi} and \eqref{con:gamma_i} of Theorem \ref{the:newtech} are satisfied.
\qed 

\begin{remark}\label{rem}
We will now show that by taking $\Phi(x)=\log(x+1)$ and special choices of the constants, the main results of \cite{DehDur11} and \cite{DehDurVol09} can be obtained as corollaries of Theorem \ref{pro:newtech_hoelder}. Note that this $\Phi$-function arises in \eqref{con:2p-moment} when the process is multiple mixing with exponential decay (see \cite{DehDur11}).
For the above choice of $\Phi$, \eqref{con:momentbound'} can be simplified to
\begin{align*} 
	w_{\F} (y) \leq \beta |\log(y)|^{-{\gamma}} 
\end{align*}
for some $\beta>0$.
\begin{enumerate}[1)]
\item If we choose $\kappa_i=2p-i$, 
we get
\begin{align*}
\max_{i=1,\ldots,p}	\frac{2p-i}{\frac{i}{r} + 2(p-i) -d} = \frac{rp}{p-rd},
\end{align*}
and thus condition \eqref{con:kappa_i} in Theorem \ref{pro:newtech_hoelder} is equivalent to 
$ \gamma > \frac{rp}{p-rd}$.
In this way, we obtain Theorem 3 in \cite[p.1081]{DehDur11}.
\item If we choose $\alpha,d,r=1$ and $p=2$, condition \eqref{con:kappa_i} in Theorem \ref{pro:newtech_hoelder} reduces to 
$ \gamma > \max\{ \frac{\kappa_1}{2}, \kappa_2 \}.$
Thus we obtain Theorem 1 in \cite[p.3702]{DehDurVol09}.
\end{enumerate}
This shows that Theorem \ref{pro:newtech_hoelder} is a far-reaching generalization of the result obtained in \cite{DehDur11} and \cite{DehDurVol09}.
\end{remark}

We finally obtain a version of Theorem \ref{the:multimix} for the special choice $\mathcal{G}=\mathcal{H}_{\alpha}$.
\begin{corollary}\label{coro:multimix}
If conditions \eqref{2-i} and \eqref{2-ii} of Theorem \ref{the:multimix} hold with $\mathcal{G}=\mathcal{H}_{\alpha}$ and if
$\F$ is $\theta$-H\"older
for some $\theta$ such that 
\begin{align}
\frac{\theta}{\alpha}> \frac{rp}{p-rd}, \label{gamma} 
 \end{align}
then the empirical CLT holds.
\end{corollary}

\proof
We will apply Theorem \ref{pro:newtech_hoelder} and thus have to verify assumption \eqref{3-ii}.
Choosing  $\kappa_i=i$ and $\Phi(x)=x$, we obtain the moment bound \eqref{con:2p-moment} by Proposition \ref{pro:multimix->2p-moment},
while \eqref{eq:phi_hoelder} is obviously satisfied.
Choosing moreover $\gamma=\frac{\theta}{\alpha}$, \eqref{con:kappa_i} follows from \eqref{gamma}. Finally 
\eqref{con:momentbound'} is a consequence of the H\"older assumption on $\F$.
\qed

\subsection{Causal functions of independent and identically distributed processes}

One example of processes that feature the multiple mixing property \eqref{con:multmix_pol} and that can be treated by our methods is the class of causal functions of i.i.d.\ processes, which are defined as follows:

\begin{definition}[Causal function]
Let $(\xi_j)_{j\in\Z}$ be an independent identically distributed process with values in a Banach space $\mathcal{X}$. We call $(X_i)_{i\in\No}$ a causal function of $(\xi_j)_{j\in\Z}$ if there is a measurable function
$G: \mathcal{X}^{\No} \ra \R^d$ such that each $X_i$ is of the form%
\begin{align*}	
	X_i:= G((\xi_{i-j})_{j\in\No}).
\end{align*}
\end{definition}

Let us now introduce a measure of the dependence structure of a causal function of an i.i.d.\ process $(\xi_j)_{j\in\Z}$. Set
\begin{align*}	
	\dot{X}_i&:= G(\xi_{i},\xi_{i-1},\ldots,\xi_{1},\xi'_{0},\xi'_{-1},\ldots),
\end{align*}
where $(\xi'_j)_{j\in\Z}$ is an independent copy of $(\xi_j)_{j\in\Z}$, i.e. $(\xi_j)_{j\in\Z}$ and $(\xi'_j)_{j\in\Z}$ are identically distributed and both processes are independent from each other. We can now define for $i\in\N$ and $s\ge 1$,
\begin{align}
	\delta_{i,s} = \|X_i - \dot{X}_i \|_s := \E \Bigl( \| X_i - \dot{X}_i \|_{\R^d}^s \Bigr)^{\frac{1}{s}}. \label{delta}
\end{align}
\change This physical dependence measure was introduced by Wu \cite{Wu05}.

\begin{proposition} \label{pro:causal->multimix}
Let $(X_i)_{i\in\No}$ be an $\R^d$-valued causal function of an i.i.d.\ process, $\alpha\in (0,1]$ and $r\in[1,\infty)$, $s\in(1,\infty]$ with $\frac{1}{r}+\frac{1}{s}=1$. Then
$$(X_i)_{i\in\No}\in\MMP{\Theta}{r}{\mathcal{H}_{\alpha}},$$ 
where $\Theta(i)=(\delta_{i,s})^{\alpha}$.
\change As a consequence of Proposition \ref{pro:multimix->2p-moment}, if 
\begin{align}
	\sum_{i=1}^{\infty} i^{2p-2} (\delta_{i,s})^{\alpha} < \infty \label{eq:wu-sum}
\end{align}
for some $p>rd$, then
the moment bound \eqref{con:2p-moment} holds for all $f\in\mathcal{H}_{\alpha}$ such that $\|f\|_{\infty}\leq1$, with $p, r$ as above 
and $\Phi_i(x)=x^i$.
\end{proposition}

The second crucial point in the assumptions of Theorem \ref{the:newtech} is the CLT \eqref{con:clt}. 
The following proposition gives a criterion for \eqref{con:clt} in the situation of a causal function of an i.i.d.\ process.

\begin{proposition}\label{pro:clt_causal}
If $(X_i)_{i\in\No}$ is an $\R^d$-valued causal function of an i.i.d.\ process and satisfies
\begin{align}
	\sum_{i=1}^{\infty} (\delta_{i,s})^{\alpha} < \infty \label{con:causal},
\end{align}
for some $\alpha\in(0,1]$, $s\in[1,\infty]$, then the CLT \eqref{con:clt} under $\mathcal{H}_{\alpha}$ holds
with $\sigma_f^2=\E(f(X_0)^2)+2\sum_{i=1}^{\infty} \E\bigl(f(X_0) f(X_i)\bigr)$.
\end{proposition}

The proofs of the two preceding propositions are given in Section \ref{sec:new_technique_applications_proof_lemma}.
As a direct application of previous results, we obtain the following one:

\begin{theorem}\label{cor:newtech_causal}
Let $(X_i)_{i\in\No}$ be an $\R^d$-valued causal function of an i.i.d.\ sequence. Assume that:
\begin{enumerate}[(i)]
\item \label{4-i}
The distribution function $\F$ of $X_0$ is $\theta$-H\"older for some $\theta\in(0,1]$.
\item \label{4-ii}
There are some $r\in[1,\infty)$, $s\in(1,\infty]$ satisfying $\frac{1}{r}+\frac{1}{s}=1$, an integer $p>rd$ and a positive constant $\alpha\in(0,1]$ satisfying \eqref{gamma} and \eqref{eq:wu-sum}.
\end{enumerate}
Then there is a centered Gaussian process $(W(t))_{t\in[-\infty,\infty]^d}$ with almost surely continuous sample paths such that $U_n \dconv W$, in the space $\D([-\infty,\infty]^d)$.
\end{theorem}

\proof
\change We will apply Corollary \ref{coro:multimix}. By Proposition \ref{pro:clt_causal} the CLT \eqref{con:clt} holds under $\mathcal{H}_{\alpha}$. Proposition \ref{pro:causal->multimix} shows that $(X_i)_{i\in\No}\in\MMP{\Theta}{r}{\mathcal{H}_{\alpha}}$, where $\Theta(i)=(\delta_{i,s})^\alpha$ and thus $\sum_{i=0}^{\infty} i^{2p-2}\Theta(i)<\infty$.
\qed

\begin{example}[Time delay vectors]
 Let $(X_i)_{i\in\No}$ be a real-valued causal function of an i.i.d.\ process. We define the time delay vector process $(Y_i)_{i\in\No}$ of dimension $d\ge 1$, by
\[
 Y_i=(X_i,\ldots,X_{i+d-1}), \quad i\in\No.
\]
If the scalar process $(X_i)_{i\in\No}$ satisfies \eqref{4-i} and \eqref{4-ii} of Theorem \ref{cor:newtech_causal}, then the Empirical Central Limit Theorem holds for the process $(Y_i)_{i\in\No}$.
\proof
Assume that $(X_i)_{i\in\No}$ satisfies \eqref{4-i} and \eqref{4-ii} of Theorem \ref{cor:newtech_causal} and let us check that the process $(Y_i)_{i\in\No}$ also satisfies these assumptions. Denote by $F_X$ the distribution function of $X_0$ and by $F_Y$ the multidimensional distribution function of $Y_0$. The marginals of $F_Y$ are all $F_X$ and then $w_{F_Y}\le d w_{F_X}$. Thus $F_Y$ is $\theta$-H\"older.
Denote by $\delta_{i,s}(X)$ and $\delta_{i,s}(Y)$ the coefficients introduced in \eqref{delta} relative respectively to $(X_i)_{i\in\No}$ and $(Y_i)_{i\in\No}$. We can see that there exists a constant $C>0$ such that for all $i\in\No$,
\[
 \delta_{i,s}(Y)\le C ( \delta_{i,s}(X)+\ldots + \delta_{i+d-1,s}(X) ).
\]
Thus, we infer that $(Y_i)_{i\in\No}$ satisfies \change\eqref{eq:wu-sum} with the same constant $\alpha$ as for $(X_i)_{i\in\No}$.
\qed
\end{example}

\begin{example}[Linear processes]
Let $(X_i)_{i\in\No}$ be a causal linear function given by 
\begin{align*}
X_i := \sum_{j\in\No} a_j \xi_{i-j},
\end{align*}
where $(\xi_j)_{j\in\Z}$ is an i.i.d.\ $\mathcal{X}$-valued process and $(a_j)_{j\in\No}$ is a family of linear operators from $\mathcal{X}$ to $\R^d$. We denote the norm of such operators by 
$$|a|=\sup \{\|a(x)\| : x\in\mathcal{X}, |x|_\mathcal{X}\le 1\}.$$
If $\|\xi_0\|_s<\infty$ for some $s>1$, if the distribution function $\F$ of $X_0$ is $\theta$-H\"older and if
\change\begin{align*}
 \sum_{j=i}^{\infty} |a_j| = \mathcal{O}(i^{-b}) 
\;\mbox{ with }
b >  \min_{p\in\No, p>rd} \frac{r}{\theta} \frac{(2p-1)p}{p-rd}
\end{align*}
for $r=\frac{s}{s-1}$, then the Empirical Central Limit Theorem holds.
\end{example}

\proof
\change Let $(\xi'_j)_{j\in\Z}$ be an independent copy of $(\xi_j)_{j\in\Z}$ and $p$ an integer which realizes the minimum in the condition on $b$. By assumption, there is an $\epsilon>0$ such that
$b>(1+\epsilon) \frac{(2p-1)rp}{\theta (p -rd)}$.
We can choose $\alpha=\theta(1+\epsilon)^{-1}\frac{p-rd}{rp}$, ensuring that \eqref{gamma} is satisfied. 
Since 
\begin{align*}
 \delta_{i,s}=\biggl\|\sum_{j=i}^\infty a_j(\xi_{i-j}-\xi'_{i-j})\biggr\|_s
             \le \|\xi_0-\xi_0'\|_s \sum_{j=i}^\infty |a_j|,
\end{align*}
we have
\begin{align*}
	i^{2p-2} (\delta_{i,s})^{\alpha} 
	\leq (2\|\xi_0\|_s)^{\alpha} i^{2p-2} \Bigl(\sum_{j=i}^{\infty} |a_j|\Bigr)^{\alpha}
	= \mathcal{O}(i^{2p-2-\alpha b}),
\end{align*}
where $2p-2-\alpha b < 2p-2 -(2p-1)= -1$, thus \eqref{eq:wu-sum} holds and Theorem \ref{cor:newtech_causal} applies. 
\qed


The example of causal linear processes has already been studied by several authors. Dedecker and Prieur \cite{DedPri07} could allow lower rates of convergence for  
$\sum_{j=i}^{\infty} |a_j|$ but required that $X_0$ has a bounded density. Wu \cite{Wu08} also assumes that the underlying i.i.d.\ process has a density. In the present paper, no assumption is made on the distribution of the $\xi_i$ except moments and the distribution of $X_0$ does not need to be absolutely continuous. In the case where $\xi_i$ are $\R$-valued and $d=1$, very weak conditions can be found in \cite{Ded10}.

\section{Proof of Theorem \ref{the:newtech}}\label{sec:new_technique_proof}
To prove the convergence of the empirical process we will use the following result obtained in \cite[Theorem 2]{DehDurVol09}. Note that this result is a modification
of a classical theorem on weak convergence; see \cite[Theorem 4.2]{Bil68}.

\begin{proposition}\label{pro:bil_erw}
Let $S$ be a separable and complete space with metric $\rho$.
If $U_n,U_n^{(m)}$ and $U^{(m)}$, $n,m\in\N$, are $S$-valued random variables such that	
	\begin{align}
		U_n^{(m)} \dconv U^{(m)} \text{ as } n\ra\infty, \text{ for all }m\in\N \text{ and }   \label{cri:bil.new.1}\\
		\lim_{m\gu}\limsup_{n\ra\infty} \P(\rho(U_n,U_n^{(m)}) > \epsilon) = 0 \text{ for all }\epsilon>0, \label{cri:bil.new.2}
	\end{align}
then there is an $S$-valued random variable $U$ such that
$U_n \dconv U$  as $n\ra\infty$ and
$U^{(m)} \dconv U$  as $m\ra\infty$.
\end{proposition}

\begin{remark}
We will apply Proposition \ref{pro:bil_erw} to the situation where $S$ is chosen as $\D([-\infty,\infty]^d)$ (equipped with the Skorokhod metric $d_0$, c.f. p.\pageref{txt:multi_cadlag}) and $U_n$ denotes the empirical process.
If one wants to prove the convergence in distribution of $U_n$ Proposition \ref{pro:bil_erw} assures us, that it is sufficient to find a process $U_n^{(m)}$ which approximates $U_n$ as $m\gu$ in the sense of \eqref{cri:bil.new.2} and to show that this process is convergent in distribution for each $m$ as $\ngu$.
\end{remark}

Following the techniques presented in \cite[p.1078 ff]{DehDur11}, we begin by introducing a partition for $[-\infty,\infty]^d$. Let ${\F}_i$ be the $i$-th marginal distribution of ${\F}$, $0=r^{(m)}_0<r^{(m)}_1<\ldots<r^{(m)}_m=1$ a partition of $[0,1]$ and set for $i\in\{1,\ldots,d\}$, $j_i\in\{0,\ldots,m\}$,
\begin{align*}
 t^{(m)}_{i,j_i} := {\F}_i\rinv(r^{(m)}_{j_i}), 
\end{align*}
where 
${\F}_i\rinv(y):= \sup\{x\in [-\infty,\infty] : {\F}_i (x) \leq y \}$.\\
Note that ${\F}_i\rinv$ are injective since ${\F}_i$ are continuous. For convenience we also define $t^{(m)}_{i,m+1} := t^{(m)}_{i,m}$.
For $j\in\{0,\ldots,m+1\}^d$ set
\begin{align*} t^{(m)}_j := (t^{(m)}_{1,j_1},\ldots,t^{(m)}_{d,j_d}) = ({\F}_1\rinv(r^{(m)}_{j_1}),\ldots,{\F}_d\rinv(r^{(m)}_{j_d})).
\end{align*}
To keep notation short, denote $(x,\ldots,x)\in[-\infty,\infty]^d$ by $\mult{x}$.

We can construct a $\mathcal{G}$-approximation of the indicator function $\1_{[-{\mult{\infty}},t^{(m)}_{j-\mult{1}}]}$ by setting
for $j\in\{1,\ldots,m\}^d$ 
\begin{align} \label{phi}
\phi^{(m)}_j:= \begin{cases} \phi_{\bigl(t^{(m)}_{j-\mult{2}},t^{(m)}_{j-\mult{1}}\bigr)},& \text{if}\ j\geq\mult{2},\\
0,& \text{if}\ j\ngeq\mult{2}, \end{cases} 
\end{align}
where $\phi_{(t^{(m)}_{j-\mult{2}},t^{(m)}_{j-\mult{1}})}\in\mathcal{G}$ satisfies \eqref{con:existence_phi} and \eqref{con:normbounding}.
Observe that $t^{(m)}_{j-\mult{2}}<t^{(m)}_{j-\mult{1}}$, since all ${\F}_i\rinv$ are injective.

To approximate the empirical distribution function we introduce
	\begin{align*} 
	{\F}_n^{(m)}(t) := \sum_{j\in\{1,\ldots,m\}^d} \Bigl( \frac{1}{n}\sum_{i=1}^n \phi^{(m)}_j (X_i)\Bigr) \1_{[t^{(m)}_{j-\mult{1}},t^{(m)}_j)}(t). 
	\end{align*}
Note that for $t$ in any fixed rectangle $[t^{(m)}_{j-\mult{1}},t^{(m)}_j)$ we have the simple form $${\F}_n^{(m)}(t) = \frac{1}{n}\sum_{i=1}^n \phi^{(m)}_j (X_i).$$ By the definition of the $\phi^{(m)}_j$ it is easy to see that therefore
\[ {\F}_n(t^{(m)}_{j-\mult{2}}) \leq {\F}_n^{(m)}(t) \leq {\F}_n(t^{(m)}_{j-\mult{1}})\ \ \forall t\in[t^{(m)}_{j-\mult{1}},t^{(m)}_j).\] 
Thus it is natural to approximate (as $m\gu$) $U_n$ by
\begin{align*}
	U_n^{(m)}&
	:= \Bigl(\sqrt{n}\bigl( {\F}_n^{(m)}(t)-{\F}^{(m)}(t)\bigr)\Bigr)_{t\in[-\infty,\infty]^d},
\end{align*}
where
\begin{align*}	
	{\F}^{(m)}(t) &:=\E\bigl({\F}_n^{(m)}(t)\bigr) 
	= \sum_{j\in\{1,\ldots,m\}^d} \E\bigl(\phi^{(m)}_j(X_0)\bigr)\1_{[t^{(m)}_{j-\mult{1}},t^{(m)}_j)}(t).
	\end{align*}

\begin{remark}
Notice that at this point the $\phi^{(m)}_j$ and thus $U_n^{(m)}$ depend heavily on the chosen partition $r_0,\ldots,r_m$ on $[0,1]$. Therefore the notation with the superscript $m$ may be misleading at first glance, but since whenever the choice of the partition matters we will only use equidistant partitions of $[0,1]$, the partitions will in all relevant situations be uniquely defined by $m$.  
\end{remark}

The central idea to prove Theorem \ref{the:newtech} is to use Proposition \ref{pro:bil_erw}. 
Thus we need to check \eqref{cri:bil.new.1} and \eqref{cri:bil.new.2} for $S=\D([-\infty,\infty]^d)$. This is done in the next two lemmas.

\begin{lemma}\label{lem:marg_conv}
For every partition $0=r^{(m)}_0<\ldots<r^{(m)}_m=1$ of $[0,1]$, $U_n^{(m)}$ converges weakly to some centered Gaussian process $W^{(m)}\in\D([-\infty,\infty]^d)$ whose sample paths are constant on each of the rectangles $[t^{(m)}_{j-\mult{1}},t^{(m)}_{j})$, $j\in\{1,\ldots,m\}^d$. 
\end{lemma}

\proof
Since all the $U_n^{(m)}$ are constant on each of the rectangles $[t^{(m)}_{j-\mult{1}},t^{(m)}_{j})$, 
it suffices to show weak convergence of the sequence of vectors
\[
 \left(\frac{1}{\sqrt{n}}\sum_{i=1}^n\left(\phi_j^{(m)}(X_i)-\E\left(\phi_j^{(m)}(X_i)\right)\right)\right)_{j\in\{1,\ldots,m\}^d},
\]
which is a consequence of the CLT \eqref{con:clt} and the Cram\'er-Wold device.
\qed

\begin{lemma}\label{lem:approx_m}
Let $0=r^{(m)}_0<r^{(m)}_1<\ldots<r^{(m)}_m=1$ be the partition of $[0,1]$ defined by $r^{(m)}_{k}=\frac{k}{m}$. Then for every $\epsilon,\eta>0$ there is a $m_0\in\N$ such that for all $m\geq m_0$
\[
	\limsup_{\ngu} \P \Bigl( \sup_{t\in[-\infty,\infty]^d} |U_n(t)-U_n^{(m)}(t)| > \epsilon \Bigr) \leq \eta.
\]
\end{lemma}

\proof 
Let us consider $\epsilon, \eta>0$ fixed for the rest of this proof. Consider the partition $0=r^{(m)}_0<\ldots<r^{(m)}_m=1$ of $[0,1]$ defined in the statement of the lemma and set 
$h=\frac{1}{m}$.
For each $k\in\N$, consider the refined partition 
\[ r^{(m)}_{q -1} = s_{q,0}^{(k)} < s_{q,1}^{(k)} < \ldots < s_{q,2^k}^{(k)} = r^{(m)}_{q} \]
of $[r_{q-1}^{(m)},r_{q}^{(m)}]$, where 
\[ s_{q,\ell}^{(k)} := r_{q-1}^{(m)} + \ell \cdot \frac{h}{2^k}, \]
$\ell\in\{0,\ldots,2^k\}$ and $q\in\{0,\ldots,m\}$. Setting for $i\in\{1,\ldots,d\}$, $j_i\in\{1,\ldots,m\}$, $l_i\in\{0,\ldots,2^k\}$
\[ s_{i,j_i,l_i}^{(k)} = {\F}_i\rinv (s_{j_i,l_i}^{(k)}), \]
we obtain partitions 
\[ t^{(m)}_{i,j_i-1} = s_{i,j_i,0}^{(k)} < s_{i,j_i,1}^{(k)} < \ldots < s_{i,j_i,2^k}^{(k)} = t^{(m)}_{i,j_i} \]
of $[t^{(m)}_{i,j_i-1},t^{(m)}_{i,j_i}]$. To simplify the notation in the following calculations we set 
\[
s_{i,j_i,-1}^{(k)}:=s_{i,j_i-1,2^k-1}^{(k)},\mbox{ for } j_i>1,\;\mbox{ and  }\; s_{i,j_i,2^k+1}^{(k)}:=s_{i,j_i+1,1}^{(k)}, \mbox{ for } j_i<m.
\]

Let us now focus on a fixed rectangle $[t^{(m)}_{j-\mult{1}},t^{(m)}_j)$ for some $j=(j_1,\ldots,j_d)\in\{1,\ldots,m\}^d$. Our aim is to construct a chain to link the point $t^{(m)}_{j-\mult{1}}$ to some arbitrary point $t\in[t^{(m)}_{j-\mult{1}},t^{(m)}_j)$. Therefore we set
\[ l_{i,j_i}(k,t)=\max\bigl\{\ell\in\{0,\ldots,2^k\} : s_{i,j_i,\ell}^{(k)} \leq t_i\bigr\}\in\{0,\ldots,2^k-1\}.\]
Since we consider $j$ to be fixed, we may drop the index $j$ in order to simplify further notation. More precisely, we set
\begin{align*}
s_l^{(k)}:= (s_{1,j_1,l_1}^{(k)},\ldots,s_{d,j_d,l_d}^{(k)}),\; \mbox{ and }\;
l(k,t):= (l_{1,j_1}(k,t),\ldots,l_{d,j_d}(k,t)).
\end{align*}
In this way for any $k\in\N$, we obtain an ($[-\infty,\infty]^d$- valued) chain
\[ 
	t^{(m)}_{j-\mult{1}} = s_{l(0,t)}^{(0)} \leq s_{l(1,t)}^{(1)} \leq \ldots \leq s_{l(k,t)}^{(k)} \leq t \leq s_{l(k,t)+\mult{1}}^{(k)}. 
\]

Now set 
\[ \psi_{\mult{0}}^{(0)} := \phi_{(t^{(m)}_{j-\mult{1}},t^{(m)}_j)}\]
and choose for every $k\in\N$ and $l\in\{0,\ldots,2^k+1\}^d$ a function $\psi_l^{(k)}\in\mathcal{G}$ such that%
\arpion{the reference to the indices $j$ and $m$ is omitted, since these are considered to be fixed.}
\begin{align}
	\psi_l^{(k)}=
	\begin{cases}
		0,	& \text{ if } \exists i\in\{1,\ldots,d\}: j_i=1 \text{ and } l_i=0\\
		1,	& \text{ if } (\exists i\in\{1,\ldots,d\}: j_i=m \text{ and } l_i=2^k+1)\\
				& \text{ and } (\nexists i\in\{1,\ldots,d\}: j_i=1 \text{ and } l_i=0)\\
		\phi_{(s_{l-\mult{1}}^{(k)},s_{l}^{(k)})} & \text{ else},
	\label{val:psi_l^k}
	\end{cases}
\end{align}
where $\phi_{(s_{l-\mult{1}}^{(k)},s_{l}^{(k)})}$ satisfies \eqref{con:existence_phi} and \eqref{con:normbounding}.

By this definition we have for every $t\in[-\infty,\infty]^d$ and $l\in\{0,\ldots,2^k\}^d$ the following inequalities: 
\begin{align}
	\1_{[-{\mult{\infty}},s_{l-\mult{1}}]} \leq \psi_l^{(k)} \leq \1_{[-{\mult{\infty}},s_l]} \label{eq:temp005}
\end{align}
\begin{align}
	\phi^{(m)}_j \leq \psi_{l(1,t)}^{(1)} \leq \ldots \leq \psi_{l(k,t)}^{(k)} \leq \1_{[-{\mult{\infty}},t]} \leq \psi_{l(k,t)+\mult{2}}^{(k)}.
	\label{eq:temp006}
\end{align}
Using inequality \eqref{eq:temp006} we obtain for $t\in[t^{(m)}_{j-\mult{1}},t^{(m)}_j)$ and $K\in\N$, the telescopic-sum representation
\begin{align}
	&\frac{1}{n}\sum_{i=1}^n  \Bigl( \1_{[-{\mult{\infty}},t]}(X_i) -{\F}_n^{(m)}(t) \Bigr) \notag\\
\CUT{&=\frac{1}{n}\sum_{i=1}^n \1_{[-{\mult{\infty}},t]}(X_i) - \frac{1}{n}\sum_{i=1}^n \phi^{(m)}_j (X_i) \notag\\}
	&= \sum_{k=1}^{K} \frac{1}{n}\!\sum_{i=1}^{n} \Bigl( \psi_{l(k,t)}^{(k)}(X_i) - \psi_{l(k-1,t)}^{(k-1)}(X_i) \Bigr) 
	+	\frac{1}{n}\!\sum_{i=1}^{n} \Bigl( \1_{[-{\mult{\infty}},t]}(X_i) - \psi_{l(K,t)}^{(K)}(X_i) \Bigr).
\label{eq:temp002}
\end{align}
Let us now consider
\begin{align*}
	U_n(t)-U_n^{(m)}(t) 
	= {\sqrt{n}} \Bigl(\frac{1}{n}\sum_{i=1}^n \1_{[-{\mult{\infty}},t]}(X_i) - {\F}(t) \Bigr)  
	-\sqrt{n}\Bigl( {\F}_n^{(m)}(t)-{\F}^{(m)}(t)\Bigr).
\end{align*}
Equation \eqref{eq:temp002} yields
\begin{align}
	&U_n(t)-U_n^{(m)}(t) \notag\\
	&= \sum_{k=1}^{K} \frac{1}{\sqrt{n}} \sum_{i=1}^{n} \bigl(\psi_{l(k,t)}^{(k)}(X_i) - \E\psi_{l(k,t)}^{(k)}(X_0)\bigr) 
	- \bigl(\psi_{l(k-1,t)}^{(k-1)}(X_i)-\E\psi_{l(k-1,t)}^{(k-1)}(X_0)\bigr) \notag \\
	&\hspace{3ex} +	\frac{1}{\sqrt{n}} \sum_{i=1}^{n}  \bigl(\1_{[-{\mult{\infty}},t]}(X_i)- {\F}(t) \bigr) - \bigl(\psi_{l(K,t)}^{(K)}(X_i)-\E\psi_{l(K,t)}^{(K)}(X_0)\bigr). \label{eq:temp007}
\end{align}
Applying the inequalities in \eqref{eq:temp006}, we gain the following upper bounds for the last sum on the right-hand side of the above inequality:
\begin{align}
	&\frac{1}{\sqrt{n}} \sum_{i=1}^{n} \Bigl( \bigl(\1_{[-{\mult{\infty}},t]}(X_i)- {\F}(t) \bigr) - \bigl(\psi_{l(K,t)}^{(K)}(X_i)-\E\psi_{l(K,t)}^{(K)}(X_0)\bigr)\Bigr) \notag\\
\CUT{  &\geq -\sqrt{n} \Bigl({\F}(t)-\E\psi_{l(K,t)}^{(K)}(X_0)\Bigl) \notag\\}
	&\geq -\sqrt{n} \Bigl(\E\psi_{l(K,t)+\mult{2}}^{(K)}(X_0)-\E\psi_{l(K,t)}^{(K)}(X_0)\Bigl)\label{eq:temp003}
\end{align}	
and
\begin{align}	
	&\frac{1}{\sqrt{n}} \sum_{i=1}^{n} \Bigl( \bigl(\1_{[-{\mult{\infty}},t]}(X_i)- {\F}(t) \bigr) - \bigl(\psi_{l(K,t)}^{(K)}(X_i)-\E\psi_{l(K,t)}^{(K)}(X_0)\bigr)\Bigr) \notag\\
\CUT{	
	&\leq \frac{1}{\sqrt{n}} \sum_{i=1}^{n} \Bigl( \bigl(\psi_{l(K,t)+\mult{2}}^{(K)}(X_i) - \E\psi_{l(K,t)+\mult{2}}^{(K)}(X_0)\bigr)- \bigl(\psi_{l(K,t)}^{(K)}(X_i)-\E\psi_{l(K,t)}^{(K)}(X_0)\bigr) \notag\\
	&\hspace{3ex} + \sqrt{n} \Bigl(\E\psi_{l(K,t)+\mult{2}}^{(K)}(X_0) - {\F}(t)\Bigl)\notag\\
	}
	&\leq \frac{1}{\sqrt{n}} \sum_{i=1}^{n} \Bigl( \bigl(\psi_{l(K,t)+\mult{2}}^{(K)}(X_i) - \E\psi_{l(K,t)+\mult{2}}^{(K)}(X_0)\bigr)- \bigl(\psi_{l(K,t)}^{(K)}(X_i)-\E\psi_{l(K,t)}^{(K)}(X_0)\bigr) \notag\\
	&\hspace{3ex} + \sqrt{n} \Bigl(\E\psi_{l(K,t)+\mult{2}}^{(K)}(X_0) - \E\psi_{l(K,t)}^{(K)}(X_0)\Bigl). \label{eq:temp004}
\end{align} 

For convenience, let $s_{i,m,2^k+1}^{(k)}:=s_{i,m,2^k}^{(k)}$. By equation \eqref{eq:temp005} and the continuity\arpion{note that for continuous ${\F}$, we have ${\F}\circ {\F}\rinv(x)=x$ for all $x\in[0,1]$.} of ${\F}$, we obtain 
\begin{align*} 
	&{\sqrt{n}}\Bigl|\E\psi_{l(K,t)+\mult{2}}^{(K)}(X_0) - \E\psi_{l(K,t)}^{(K)}(X_0)\Bigr| \\
	&\leq{\sqrt{n}} \Bigl|\E\1_{[-{\mult{\infty}},s_{l(K,t)+\mult{2}})}(X_0)-\E\1_{[-{\mult{\infty}},s_{l(K,t)-\mult{1}})}(X_0)\Bigr| \\
\CUT{	&= {\sqrt{n}} \Bigl({\F}(s_{l(K,t)+\mult{2}})-{\F}(s_{l(K,t)-\mult{1}})\Bigr)\\
	}
	&\leq{\sqrt{n}} \Bigl( d\max_{i=1,\ldots,d} \bigl\{ {\F}_i(s_{i,j_i,l_{i,j_i}(K,t)+\mult{2}})-{\F}_i(s_{i,j_i,l_{i,j_i}(K,t)-\mult{1}}) \bigr\} \Bigr)
	= \frac{3d\sqrt{n}h}{2^K}, 
\end{align*}
 and thus, if we choose
\begin{align}
	K=K_n:= \biggl\lfloor \log_2 \biggl( \frac{2^4 d}{\epsilon} \sqrt{n}h\biggr)\biggl\rfloor, \label{val:K_n}
\end{align}
we obtain 
\begin{align} \Bigl|\E\psi_{l(K,t)+\mult{2}}^{(K)}(X_0) - \E\psi_{l(K,t)}^{(K)}(X_0)\Bigr|< \frac{\epsilon}{2}. \label{eq:temp004.1}\end{align}

In summary, using \eqref{eq:temp003}, \eqref{eq:temp004} and \eqref{eq:temp004.1} in equation \eqref{eq:temp007} yields, for all $n\in\N$,
\begin{align}
	\bigl|U_n(t)-U_n^{(m)}(t)\bigr| 
	&< \biggl| \sum_{k=1}^{K_n} \frac{1}{\sqrt{n}} \sum_{i=1}^{n} \biggl\{ \Bigl(\psi_{l(k,t)}^{(k)}(X_i) - \E\psi_{l(k,t)}^{(k)}(X_0)\Bigr) \notag\\
	&\hspace{3ex} - \Bigl(\psi_{l(k-1,t)}^{(k-1)}(X_i)-\E\psi_{l(k-1,t)}^{(k-1)}(X_0)\Bigr) \biggr\} \biggr|\notag \\
	&\hspace{3ex} + \biggl| \frac{1}{\sqrt{n}} \sum_{i=1}^{n} \biggl\{ \Bigl(\psi_{l(K_n,t)+\mult{2}}^{(K_n)}(X_i) - \E\psi_{l(K_n,t)+\mult{2}}^{(K_n)}(X_0)\Bigr) 
	\notag\\
	&\hspace{3ex} -\Bigl(\psi_{l(K_n,t)}^{(K_n)}(X_i)-\E\psi_{l(K_n,t)}^{(K_n)}(X_0)\Bigr) \biggr\} \biggr| 
	+ \frac{\epsilon}{2}. \label{eq:temp008}
\end{align}

Consider the maximum of the terms in \eqref{eq:temp008} over all $t\in[t^{(m)}_{j-\mult{1}},t^{(m)}_j)$. By the definition of the $l(k,t)$ we have
\begin{align*} 
	\biggl\lfloor \frac{l(k,t)}{2} \biggr\rfloor :&= \Biggl(\biggl\lfloor \frac{l_{1,j_1}(k,t)}{2} \biggr\rfloor,\ldots,\biggl\lfloor \frac{l_{d,j_d}(k,t)}{2} \biggr\rfloor\Biggr) = l(k-1,t).
\end{align*}
We therefore obtain
\begin{align*}
	\sup_{t\in[t^{(m)}_{j-\mult{1}},t^{(m)}_j)} |U_n(t)-U_n^{(m)}(t)| 
		<& \sum_{k=1}^{K_n} \frac{1}{\sqrt{n}} \max_{l\in\{0,\ldots,2^{k}-1\}^d} \biggl| \sum_{i=1}^{n} \biggl\{ \Bigl(\psi_{l}^{(k)}(X_i) - \E\psi_{l}^{(k)}(X_0)\Bigr) 
	\notag\\
	&\hspace{3ex}- \Bigl(\psi_{\lfloor l/2\rfloor}^{(k-1)}(X_i)-\E\psi_{\lfloor l/2\rfloor}^{(k-1)}(X_0)\Bigr) \biggr\} \biggr| \notag \\
	& +\frac{1}{\sqrt{n}} \max_{l\in\{0,\ldots,2^{K_n}-1\}^d} \biggl| \sum_{i=1}^{n} \biggl\{ \Bigl(\psi_{l+\mult{2}}^{(K_n)}(X_i) - \E\psi_{l+\mult{2}}^{(K_n)}(X_0)\Bigr)
	\notag\\
	&\hspace{3ex}- \Bigl(\psi_{l}^{(K_n)}(X_i)-\E\psi_{l}^{(K_n)}(X_0)\Bigr) \biggr\} \biggr| \notag\\
	&+ \frac{\epsilon}{2}.
\end{align*}

Choose $\epsilon_k=\frac{\epsilon}{4k(k+1)}$ and note that $\sum_{k=1}^{\infty}\epsilon_k=\epsilon/4$.
An application of Markov's inequality for the $2p$-th moments combined with condition \eqref{con:2p-moment} implies 
\begin{align}
	&\P\biggl( \sup_{t\in[t^{(m)}_{j-1},t^{(m)}_j)} |U_n(t)-U_n^{(m)}(t)| \geq \epsilon \biggr) 
				\notag\\
	&\leq \sum_{k=1}^{K_n}\sum_{l\in\{0,\ldots,2^{k}-1\}^d} \P\biggl(  \frac{1}{\sqrt{n}} \biggl| \sum_{i=1}^{n} \biggl\{ \bigl(\psi_{l}^{(k)}(X_i) - \E\psi_{l}^{(k)}(X_0)\bigr) 
	\notag\\&\hspace{8ex}
	- \bigl(\psi_{\lfloor l/2\rfloor}^{(k-1)}(X_i)-\E\psi_{\lfloor l/2\rfloor}^{(k-1)}(X_0)\bigr) \biggr\} \biggr| > \epsilon_k \biggr) 
				\notag\\
	&\hspace{3ex}+\!\!\! \sum_{l\in\{0,\ldots,2^{K_n}-1\}^d} \! \P\biggl( \frac{1}{\sqrt{n}} \biggl| \sum_{i=1}^{n} \biggl\{ \bigl(\psi_{l+\mult{2}}^{(K_n)}(X_i) - \E\psi_{l+\mult{2}}^{(K_n)}(X_0)\bigr) 
	\notag\\&\hspace{8ex}
	- \bigl(\psi_{l}^{(K_n)}(X_i)-\E\psi_{l}^{(K_n)}(X_0)\bigr) \biggr\} \biggr| > \frac{\epsilon}{4} \biggr) 
				\allowdisplaybreaks[1]\notag\\
	\CUTT{&\leq \sum_{k=1}^{K_n}\sum_{l\in\{0,\ldots,2^{k}-1\}^d} \frac{1}{\epsilon_k^{2p}n^p} \E\biggl(  \Bigl| \sum_{i=1}^{n} \bigl(\psi_{l}^{(k)}(X_i) - \E\psi_{l}^{(k)}(X_0)\bigr) 
	\notag\\&\hspace{8ex}
	- \bigl(\psi_{\lfloor l/2\rfloor}^{(k-1)}(X_i)-\E\psi_{\lfloor l/2\rfloor}^{(k-1)}(X_0)\bigr) \Bigr|^{2p} \biggr) 
				\notag\\
	&\hspace{3ex}+\!\!\! \sum_{l\in\{0,\ldots,2^{K_n}-1\}^d} \frac{4^{2p}}{n^p\epsilon^{2p}} \E\biggl( \Bigl| \sum_{i=1}^{n} \bigl(\psi_{l+\mult{2}}^{(K_n)}(X_i) - \E\psi_{l+\mult{2}}^{(K_n)}(X_0)\bigr)
	\notag\\&\hspace{8ex}
	- \bigl(\psi_{l}^{(K_n)}(X_i)-\E\psi_{l}^{(K_n)}(X_0)\bigr)\Bigr|^{2p} \biggr) 
				\allowdisplaybreaks[1]\notag\\}
	&\leq 2C \Biggl\{ \sum_{k=1}^{K_n} \sum_{l\in\{0,\ldots,2^k-1\}^d} \frac{1}{\epsilon_k^{2p} n^p} 
	\notag\\&\hspace{8ex}
	\cdot\sum_{i=1}^p n^i \: \| \psi_{l}^{(k)}(X_0)-\psi_{\lgauss{l/2}}^{(k-1)}(X_0) \|_r^i \: \Phi_i \Bigl(2\|\psi_{l}^{(k)}-\psi_{\lgauss{l/2}}^{(k-1)}\|_{\mathcal{G}}\Bigr) 
				\notag\\
	&\hspace{3ex}+  \sum_{l\in\{0,\ldots,2^{K_n}-1\}^d} \frac{4^{2p}}{\epsilon^{2p}n^p} 
	\notag\\&\hspace{8ex}
	\cdot\sum_{i=1}^p n^i \: \| \psi_{l+\mult{2}}^{({K_n})}(X_0)-\psi_{l}^{({K_n})}(X_0) \|_r^i \: \Phi_i \Bigl(2\|\psi_{l+\mult{2}}^{({K_n})}-\psi_{l}^{({K_n})}\|_{\mathcal{G}}\Bigr)  \Biggr\}.
				\allowdisplaybreaks[1]\label{eq:temp010}
\end{align}

The critical part in these terms is the argument in the functions $\Phi_i$. We therefore collect the necessary auxiliary calculations in the following lemma:
\begin{lemma}\label{lem:temp001} For all $l\in\{0,\ldots,2^k-1\}^d$, $k,n,r\in\N$
\begin{align}
	\| \psi_{l}^{(k)} (X_0) - \psi_{\lgauss{l/2}}^{(k-1)}(X_0) \|_r &\leq \Bigl( \frac{3dh}{2^k} \Bigr)^{\frac{1}{r}}, 
\notag\\
	\| \psi_{l+\mult{2}}^{(K_n)} (X_0) - \psi_{l}^{(K_n)}(X_0) \|_r &\leq \Bigl( \frac{3dh}{2^{K_n}} \Bigr)^{\frac{1}{r}},
\notag\\
	\|\psi_{l}^{(k)}\|_{\mathcal{G}} &\leq \max \biggl\{ \Psi\Bigl(\frac{2^{k}}{h}\Bigr), \|1\|_{\mathcal{G}} \biggr\} \notag
\end{align}
\end{lemma}
\proof  
By \eqref{eq:temp005} and the continuity of the ${\F}_i$,
\begin{align*}	
	\|\psi_{l}^{(k)} (X_0) - \psi_{\lgauss{l/2}}^{(k-1)}(X_0) \|_r 
	&\leq \|\1_{[-{\mult{\infty}},s_{l}^{(k)}]} (X_0) - \1_{[-{\mult{\infty}},s_{\lgauss{l/2}-\mult{1}}^{(k-1)}]}(X_0) \|_r \\
\CUT{	&\leq \biggl( \sum_{i=1}^d \bigl| {\F}_i(s_{i,j_i,l_i}^{(k)})-{\F}_i(s_{i,j_i,\lfloor l_i/2\rfloor-1}^{(k-1)}) \bigr| \biggr)^{\frac{1}{r}} \\
}
	& \leq \biggl(d \max_{i=1,\ldots,d} \bigl( {\F}_i(s_{i,j_i,l_i}^{(k)})- {\F}_i(s_{i,j_i,l_i-3}^{(k)}) \bigr)\biggr)^{\frac{1}{r}} \\
	&\leq \Bigl( \frac{3dh}{2^k} \Bigr)^{\frac{1}{r}}.
\end{align*}

The second inequality can be proven in a similar way.

In the first two cases of the definition \eqref{val:psi_l^k}, $\psi_l^{(k)}$ is a constant function taking either the value zero or one for each argument.
In this cases the last inequality of the lemma is trivially satisfied by the conditions on $\|\cdot\|_{\mathcal{G}}$. 
Else $\psi_l^{(k)}$ has a representation $\phi_{(s_{l-\mult{1}}^{(k)},s_{l}^{(k)})}$, where
\begin{align*}
s_{l}^{(k)} 
&=\Bigl(s_{1,j_1,l_1}^{(k)},\ldots,s_{d,j_d,l_d}^{(k)}\Bigr)
=\Bigl({\F}_1\rinv\bigl(s_{j_1,l_1}^{(k)}\bigr),\ldots,{\F}_d\rinv\bigl(s_{j_d,l_d}^{(k)}\bigr)\Bigr), \\
s_{l-\mult{1}}^{(k)}
&=\Bigl(s_{1,j_1,l_1-1}^{(k)},\ldots,s_{d,j_d,l_d-1}^{(k)}\Bigr)\\
&=\Bigl({\F}_1\rinv\bigl(s_{j_1,l_1}^{(k)}-h2^{-k}\bigr),\ldots,{\F}_d\rinv\bigl(s_{j_d,l_d}^{(k)}-h2^{-k}\bigr)\Bigr)
\end{align*}
and hence, for every $i\in\{1,\ldots,d\}$,
\begin{align}
	s_{i,j_i,l_i}^{(k)}-s_{i,j_i,l_i-1}^{(k)}
	\in &\{ \delta >0 : \exists t\in\R,\ |{\F}_i(t)-{\F}_i(t-\delta)|\geq h2^{-k} \} \notag\\
	&\subset \{ \delta>0 : w_{{\F}_i}(\delta) \geq h2^{-k}\} \label{eq:temp034}
\end{align}
To see this, set $\delta={\F}_i\rinv(s_{l_i}^{(k)})-{\F}_i\rinv(s_{l_i-1}^{(k)})>0$, $t={\F}_i\rinv(s_{l_i}^{(k)})$ and recall that the $\F_i\rinv$ are injective.
Now condition \eqref{con:normbounding} yields
$$\| \phi_l^{(k)} \|_{\mathcal{G}}  \leq \Psi\biggl( \frac{1}{\underset{i=1,\ldots,d}{\min} w_{\F_i} (s_{i,j_i,l_i}^{(k)}-s_{i,j_i,l_i-1}^{(k)})} \biggr) \leq \Psi\Bigl(\frac{2^k}{h}\Bigr),$$
since $\min_{i=1,\ldots,d} w_{\F_i} \bigl(s_{j_i,l_i}^{(k)}-s_{1,j_i,l_i-1}^{(k)}\bigr)\geq h2^{-k}$ by \eqref{eq:temp034}.
\qed

\medskip

An application of Lemma \ref{lem:temp001} to \eqref{eq:temp010} yields
\begin{align}
	&\P\biggl( \sup_{t\in[t^{(m)}_{j-1},t^{(m)}_j)} |U_n(t)-U_n^{(m)}(t)| \geq \epsilon \biggr) 
				\allowdisplaybreaks[1]\notag\\
	&\leq 2C \Biggl\{ \sum_{k=1}^{K_n} \sum_{i=1}^p \frac{2^{dk} n^{-(p-i)}}{\epsilon_k^{2p}}   \Bigl(\frac{3dh}{2^k}\Bigr)^{\frac{i}{r}} {\Phi_i} \Bigl(2{\Psi}\Bigl(\frac{2^k}{h}\Bigr)\Bigr) 
	\notag\\
	&\hspace{3ex}+  \sum_{i=1}^p \frac{2^{d{K_n}} n^{-(p-i)}}{(\frac{\epsilon}{4})^{2p}}  \Bigl(\frac{3dh}{2^{K_n}}\Bigr)^{\frac{i}{r}} {\Phi_i} \Bigl(2{\Psi}\Bigl(\frac{2^{K_n}}{h}\Bigr)\Bigr) \Biggr\}
	\allowdisplaybreaks[1]\notag\\
\CUT{	&\leq 2 C 
\sum_{i=1}^p \Biggl\{ (3d)^{\frac{i}{r}} n^{-(p-i)} \sum_{k=1}^{K_n} \frac{2^{dk}}{\epsilon_k^{2p}} \Bigl(\frac{h}{2^k}\Bigr)^{\frac{i}{r}} {\Phi_i} \Bigl(2{\Psi}\Bigl(\frac{2^k}{h}\Bigr)\Bigr) \Biggr\} 
	\allowdisplaybreaks[1]\notag\\
	}
	&\leq D \sum_{i=1}^p \Biggl\{ n^{-(p-i)} \sum_{k=1}^{K_n} 2^{(d-\frac{i}{r})k} k^{4p} {\Phi_i} \Bigl(2{\Psi}\Bigl(\frac{2^k}{h}\Bigr)\Bigr) h^{\frac{i}{r}} \Biggr\} 
	\allowdisplaybreaks[1]\notag\\
	&\leq D \sum_{i=1}^{p-1} \Biggl\{ n^{-(p-i)} \biggl(\frac{2^{K_n}}{h}\biggr)^{d-\frac{i}{r}} {K_n}^{4p+1} {\Phi_i} \biggl(2{\Psi}\biggl(\frac{2^{K_n}}{h}\biggr)\biggr) h^d \Biggr\}  
	\notag\\
	&\hspace{3ex} + D \Biggl\{ \sum_{k=1}^{K_n} 2^{(d-\frac{p}{r})k} k^{4p}  {\Phi_p} \Bigl(2{\Psi}\Bigl(\frac{2^k}{h}\Bigr)\Bigr) h^{\frac{p}{r}} \Biggr\} \label{eq:temp011}
\end{align}
for every $j\in\{1,\ldots,m\}^d$, where $D>0$ denotes some finite constant. In the second inequality we used that  ${\Psi}$ and ${\Phi_i}$ are nondecreasing functions and $\epsilon/4 > \epsilon_{K_n}$.

Let us first deal with the term in the last line of \eqref{eq:temp011}. By condition \eqref{con:Phi_of_Psi} we have 
\begin{align*}
	 \sum_{k=1}^{K_n} 2^{(d-\frac{p}{r})k} k^{4p} {\Phi_p} \Bigl(2{\Psi}\Bigl(\frac{2^k}{h}\Bigr)\Bigr) h^{\frac{p}{r}} 
	\leq  C'' {h}^{\frac{p}{r}-\gamma_p} \sum_{k=1}^{\infty} 2^{(\gamma_p-(\frac{p}{r}-d))k} k^{4p} 
\end{align*}
where $\gamma_p<\frac{p}{r}-d$. Hence there is a nonnegative constant $D'<\infty$ such that
\begin{align}
		\sum_{k=1}^{K_n} 2^{(d-\frac{p}{r})k} k^{4p} {\Phi_p} \Bigl(2{\Psi}\Bigl(\frac{2^k}{h}\Bigr)\Bigr) h^{\frac{p}{r}} \leq D' {h}^{\frac{p}{r}-\gamma_p} = o(h^d). \label{eq:temp012}
\end{align}

Now consider the first summand on the right-hand side of inequality \eqref{eq:temp011}. In \eqref{val:K_n} we chose $K_n= \lfloor \log_2 ({2^4 d\sqrt{n}h}/{\epsilon}) \rfloor$, hence condition \eqref{con:Phi_of_Psi} yields for any $i=1,\ldots,p-1$,
\begin{align*}
	&n^{-(p-i)} \biggl(\frac{2^{K_n}}{h}\biggr)^{d-\frac{i}{r}} {K_n}^{4p+1} {\Phi_i} \biggl(2{\Psi}\biggl(\frac{2^{K_n}}{h}\biggr)\biggr) h^d \\
	&\leq D'' \log_2^{4p+1} \Bigl( \frac{2^4 d}{\epsilon} \sqrt{n}h \Bigr) 
	\cdot (\sqrt{n})^{\gamma_i-(\frac{i}{r} + 2(p-i) - d)} h^{d}
\end{align*}
for some nonnegative constant $D''<\infty$. Since $\gamma_i< \frac{i}{r} + 2(p-i) - d$ for $i=1,\ldots,p-1$, by \eqref{con:Phi_of_Psi} we obtain for all $\eta>0$ and sufficiently large $n\in\N$,
\begin{align}
	&D \sum_{i=1}^{p-1} \Biggl\{ n^{-(p-i)} \biggl(\frac{2^{K_n}}{h}\biggr)^{d-\frac{i}{r}} {K_n}^{4p+1} {\Phi_i} \biggl(2{\Psi}\biggl(\frac{2^{K_n}}{h}\biggr)\biggr) h^d \Biggr\}  
 \leq \frac{1}{2}\eta h^{d}. \label{eq:temp013}
\end{align}

Finally, by \eqref{eq:temp011}, \eqref{eq:temp012} and \eqref{eq:temp013}, for any $\eta>0$ 
\begin{align*}
	&\limsup_{\ngu}\P\biggl( \sup_{t\in[-\infty,\infty]^d} |U_n(t)-U_n^{(m)}(t)| \geq \epsilon \biggr)\\ 
	&\leq	\limsup_{\ngu} \sum_{j\in\{1\ldots,m\}^d} \P\biggl( \sup_{t\in[t^{(m)}_{j-1},t^{(m)}_j)} |U_n(t)-U_n^{(m)}(t)| \geq \epsilon \biggr) \\	
	&\leq m^d \Bigl( o(h^d) + \frac{1}{2}\eta h^d \Bigr) = m^d \Bigl( o(m^{-d}) + \frac{1}{2}\eta m^{-d} \Bigr), 
\end{align*}
since $h=1/m$. Hence there is a $m_0\in\N$ such that 
\begin{align*}	
	\limsup_{\ngu}\P\biggl( \sup_{t\in[-\infty,\infty]^d} |U_n(t)-U_n^{(m)}(t)| \geq \epsilon \biggr) \leq \eta			
\end{align*}
for all $m\geq m_0$.
\qed

With Lemma \ref{lem:marg_conv} and Lemma \ref{lem:approx_m} established, let us finally prove Theorem \ref{the:newtech}.

\proof[Theorem \ref{the:newtech}]
By application of Proposition~\ref{pro:bil_erw} on $\D([-\infty,\infty]^d)$ equipped with the Skorokhod metric $\rho$, Lemma~\ref{lem:marg_conv} (with $r_k^{(m)}:=\frac{k}{m}$) and Lemma~\ref{lem:approx_m} show that $U_n$ converges in distribution to a process $W$ which is also the limit process of the sequence $W^{(m)}$, $m\in\N$. Since all $W^{(m)}$ are centered Gaussian processes the limit process must also be centered Gaussian. 

It remains to prove the continuity of the sample paths of $W$. At this point we already know that $U_n$ converges weakly to $W$. Therefore it is sufficient to show that for every $\epsilon, \eta>0$, there is a $\delta>0$ such that 
\begin{align}\label{eq:sufficient_condition_for_continuity}
 \limsup_{\ngu} \P \Bigl( \sup_{\|t-s\| < \delta} | U_n(t)-U_n(s)| > 3\epsilon  \Bigr) < 3 \eta.
\end{align}
The sufficiency of this condition can be proven exactly the same way as in the proof of Theorem 15.5 in \cite[p.127 f.]{Bil68}.

For all $m\in\N$, by some triangle inequality arguments we obtain
\begin{align}
	&\limsup_{\ngu} \P \Bigl( \sup_{\|t-s\| < \delta} | U_n(t)-U_n(s)| > 3\epsilon  \Bigr) \notag\\
	&\leq 2 \limsup_{\ngu} \P \Bigl( \sup_{t} | U_n(t)-U^{(m)}_n(t)| > \epsilon  \Bigr) \notag\\
	&\hspace{3ex} + \limsup_{\ngu} \P \Bigl( \sup_{\|t-s\| < \delta} | U^{(m)}_n(t)-U^{(m)}_n(s)| > \epsilon  \Bigr) \notag
\end{align}
and thus, by Lemma \ref{lem:approx_m}, there is a $m_0\in\N$ such that for all $m\geq m_0$,
\begin{align}	
	&\limsup_{\ngu} \P \Bigl( \sup_{\|t-s\| < \delta} | U_n(t)-U_n(s)| > 3\epsilon  \Bigr) \notag\\
	&\leq 2 \eta
	+ \limsup_{\ngu} \P \Bigl( \sup_{\|t-s\| < \delta} | U^{(m)}_n(t)-U^{(m)}_n(s)| > \epsilon  \Bigr). \label{eq:tmp015}
\end{align}
Now set $\delta_m := \frac{1}{2} \min_{j\in\{0,\ldots,m\}^d}\Bigl\{ \max_{i=1,\ldots,d} |t_{j_i}-t_{{j_i}-1}| \Bigr\}$
and observe that $\delta_m$ is strictly positive for any $m\in\N$, since the $\F_i\rinv$ used in the construction of the $t_j$ are strictly increasing.
Obviously for all $\delta\leq\delta_m$ and $\|t-s\|<\delta$, the points $s,t\in[-\infty,\infty]^d$ must be located in adjacent (or identical) intervals of the form $[t_j,t_{j-\mult{1}})$. 
Since the process $U_n^{(m)}$ is constant on any of the intervals $[t_j,t_{j-\mult{1}})$ and by symmetry in the arguments $s,t$ we obtain 
\[
\sup_{\|t-s\| < \delta} | U^{(m)}_n(t)-U^{(m)}_n(s)| = \max_{\substack{j\in\{0,\ldots,m\}^d \\ z\in\{0,1\}^d,\ j\geq z}} | U^{(m)}_n(t_j)-U^{(m)}_n(t_{j-z})|,
\]
thus 
\begin{align} 
	& \P \Bigl( \sup_{\|t-s\| < \delta} | U^{(m)}_n(t)-U^{(m)}_n(s)| > \epsilon  \Bigr) \notag\\
	&\leq 2^d (m+1)^d \max_{\substack{j\in\{0,\ldots,m\}^d \\ z\in\{0,1\}^d,\ j\geq z}} 
	\P \Bigl( | U^{(m)}_n(t_j)-U^{(m)}_n(t_{j-z})| > \epsilon  \Bigr). \label{eq:temp014}
\end{align}
Recall that the functions $\phi^{(m)}_{j}$ are defined in \eqref{phi}.
Analogously to the calculations in Lemma~\ref{lem:temp001}, one can show that for all $j\in\{0,\ldots,m\}^d$ and $z\in\{0,1\}^d$ such that $j\geq z$, we have
\begin{align*}
	\| \phi_{j+\mult{1}}^{(m)} (X_0) - \phi^{(m)}_{j+\mult{1}-z}(X_0) \|_r &\leq \Bigl( \frac{3d}{m} \Bigr)^{\frac{1}{r}}, 
\notag\\
	\|\phi_{j+\mult{1}}^{(m)}\|_{\mathcal{G}} &\leq \max \Bigl\{ \Psi(m), \|1\|_{\mathcal{G}} \Bigr\}.
\end{align*}
Then, by applying one after another Markov's inequality, the $2p$-th moment bounds \eqref{con:2p-moment} and the preceding inequalities, we obtain
\begin{align*}
	&\P \Bigl( | U^{(m)}_n(t_j)-U^{(m)}_n(t_{j-z})| > \epsilon  \Bigr) \\
\CUT{	
&=\P \biggl(  \biggl| \sum_{i=1}^n \biggl\{ \bigl( \phi^{(m)}_{j+\mult{1}}(X_i) -\phi^{(m)}_{j+\mult{1}-z}(X_i) \bigr) 
	- \E\bigl( \phi^{(m)}_{j+\mult{1}}(X_0) -\phi^{(m)}_{j+\mult{1}-z}(X_0) \bigr) \biggr\} \biggr|^{2p} > n^{p}\epsilon^{2p}  \biggr) \allowdisplaybreaks[1]\\
	}
	&\leq n^{-p}\epsilon^{-2p} \E\Bigl| \sum_{i=1}^n \bigl( \phi^{(m)}_{j+\mult{1}}(X_i) -\phi^{(m)}_{j+\mult{1}-z}(X_i) \bigr) - \E\bigl( \phi^{(m)}_{j+\mult{1}}(X_0) -\phi^{(m)}_{j+\mult{1}-z}(X_0) \bigr) \Bigr|^{2p}\allowdisplaybreaks[1]\\
	&\leq 2C n^{-p} \epsilon^{-2p} \sum_{i=1}^p n^{i} \|\phi^{(m)}_{j+\mult{1}}(X_0)-\phi^{(m)}_{j+\mult{1}-z}(X_0)\|_r^i \Phi_i\Bigl( 2\|\phi^{(m)}_{j+\mult{1}}-\phi^{(m)}_{j+\mult{1}-z}\|_{\mathcal{G}} \Bigr)\allowdisplaybreaks[1] \\
	&\leq 2C n^{-p} \epsilon^{-2p} \sum_{i=1}^p n^{i} \Bigl( \frac{3d}{m} \Bigr)^{\frac{i}{r}} \Phi_i\Bigl(2\Psi(m) \Bigr) \\
\CUT{	&\leq D \sum_{i=1}^p n^{-(p-i)} m^{\gamma_i-\frac{i}{r}} \\
	}
	&\leq D m^{\gamma_p-(\frac{p}{r})} + D \sum_{i=1}^{p-1} n^{-(p-i)} m^{\gamma_i-\frac{i}{r}},
\end{align*}
where $D$ is some finite constant. Therefore by \eqref{eq:temp014} there is another finite constant $D'$ such that
\begin{align*}
	& \P \Bigl( \sup_{\|t-s\| < \delta} | U^{(m)}_n(t)-U^{(m)}_n(s)| > \epsilon  \Bigr)
	\leq D' m^d \Bigl( m^{\gamma_p-\frac{p}{r}} + \sum_{i=1}^{p-1} n^{-(p-i)} m^{\gamma_i-\frac{i}{r}} \Bigr), 	
\end{align*}
and thus
\begin{align*}
	\limsup_{\ngu} & \P \Bigl( \sup_{\|t-s\| < \delta_{m}} | U^{(m)}_n(t)-U^{(m)}_n(s)| > \epsilon  \Bigr)
	\leq D' m^{\gamma_p-(\frac{p}{r}-d)} < \eta 	
\end{align*}
for sufficiently large $m\in\N$, say $m\geq m_1$. 
By \eqref{eq:tmp015} this implies that \eqref{eq:sufficient_condition_for_continuity} holds for $\delta=\delta_{\max\{m_0,m_1\}}$. 
\qed	

\begin{remark}\label{rem:class_G}
We saw in the proof, that the theorem also holds if \eqref{con:2p-moment} is only satisfied for a certain subclass of functions in $\mathcal{G}$, more precisely if \eqref{con:2p-moment} holds for all functions $f\in\mathcal{G}$ of the form
\begin{align*}
	f:= \phi_{(a,b)}- \phi_{(a',b')},
\end{align*}
where $a,b,a',b'\in[-\infty,\infty]^d$, $a'<b$, are such that

\begin{align} P(X_0\in[a',b']) \leq 2 P(X_0\in[a,b]) \leq P(X_0\in[a',b]) \leq 3 P(X_0\in[a',b']). \label{eq:temp026}\end{align}
Choosing for each $m\in\N$ an $f_m:=\phi_{(a,b)}- \phi_{(a',b')}$ such that \eqref{eq:temp026} is satisfied for $P(X_0\in[a',b'])=1/m$, it can be shown that 
\[ \| f_m(X_0) \|_r^i = \mathcal{O}\Bigl(\frac{m^{\gamma_i-\frac{i}{r}}}{\Phi_i(\|f_m\|_{\mathcal{G}})} \Bigr)\ \ \text{as}\ m\gu. \]
\end{remark}

\section{Proof of Proposition \ref{pro:multimix->2p-moment}, Proposition \ref{pro:causal->multimix} and Proposition \ref{pro:clt_causal}}\label{sec:new_technique_applications_proof}\label{sec:new_technique_applications_proof_lemma}
\subsection{Proof of Proposition \ref{pro:multimix->2p-moment}}	
By stationarity, we have 
\begin{align*}
	\Bigl| \E\Bigl( \Bigl( \sum_{i=1}^n f(X_i) \Bigr)^{p} \Bigr) \Bigr|
	&= \Bigl| \sum_{1\leq i_1,\ldots, i_p \leq n}\E \bigl( f(X_{i_1})\cdot\ldots\cdot f(X_{i_p}) \bigr)\Bigr| \notag\\
	&\leq p!n \Bigl| \sum_{\substack{0\leq i_1,\ldots, i_{p-1} \leq n-1 \\ i_{p-1}^{\ast}\leq n-1}} \E \bigl(f(X_0)f(X_{i_1^{\ast}})\cdot\ldots\cdot f(X_{i_{p-1}^{\ast}}) \bigr)\Bigr| . \notag 
\end{align*}
Using the notations $I_n(0):= \bigl|\E(f(X_0))\bigr|=0$ and
\begin{align}
	I_n(p)&:= \sum_{\substack{0\leq i_1,\ldots,i_p\leq n-1 \\ i_p^{\ast} \leq n-1}} 
	\bigl| \E \bigl( f(X_0)f(X_{i_1^{\ast}}) \cdot\ldots\cdot f(X_{i_p^{\ast}}) \bigr) \bigr|, \label{eq:temp021}
\end{align}
for $p\in\N$, we therefore have
\begin{align}
	\Bigl| \E\Bigl( \Bigl( \sum_{i=1}^n f(X_i) \Bigr)^{p} \Bigr) \Bigr|
	\leq p!n I_n(p-1).\label{eq:temp024}
\end{align}
Decomposing the sum in \eqref{eq:temp021} with respect to the highest increment of indices $i_q$, $q\in\{1,\ldots,p\}$, we receive a bound 
\begin{align*}
	I_n(p) \leq \sum_{q=1}^{p} J_n(p,q), 
\end{align*}
where
\begin{align*}
	J_n(p,q)&= \sum_{i_q=0}^{n-1} \sum_{\substack{0\leq i_1,\ldots,i_{q-1},i_{q+1},\ldots,i_p\leq i_q \\ i_p^{\ast} \leq n-1}} 
	\bigl| \E \bigl( f(X_0)f(X_{i_1^{\ast}}) \cdot\ldots\cdot f(X_{i_p^{\ast}}) \bigr) \bigr|.
\end{align*}

\begin{lemma}\label{lem:tmp_bound_Jpq} Let $(X_n)_{n\in\No}\in\MMP{\Theta}{r}{\mathcal{G}}$, then for all $p\in\N$ such that
\begin{align}
	\sum_{i=0}^{\infty} i^{p-1}\Theta(i)<\infty \label{eq:temp022}
\end{align}
and all $q\in\{1,\ldots,p\}$ there is a constant $K'$ such that 
\begin{align*}
	J_n(p,q) \leq K' \| f(X_0)\|_r \|f\|_{\mathcal{G}} + n I_n(q-1)I_n(p-q) 
\end{align*}
for all $n\in\N$ and $f\in\mathcal{G}$. 
\end{lemma}

\proof
Set
\begin{align*}
	A_{i_1,\ldots,i_p} :=& \bigl| \Cov\bigl({f(X_0)f(X_{i_1^{\ast}}) \cdot\ldots\cdot f(X_{i_{q-1}^{\ast}})}\,,\,{f(X_{i_q^{\ast}})f(X_{i_{q+1}^{\ast}}) \cdot\ldots\cdot f(X_{i_p^{\ast}})} \bigr) \bigr| \\
	B_{i_1,\ldots,i_p} :=& 
	\bigl|\E\bigl({f(X_0)f(X_{i_1^{\ast}})  \!\cdot\!\ldots\!\cdot\! f(X_{i_{q-1}^{\ast}}\!)}\bigr)\bigr|
	\!\cdot\! \bigl|\E\bigl({f(X_{0})f(X_{i_{q+1}}\!) \!\cdot\!\ldots\!\cdot\! f(X_{i_{p}^{\ast}-{i_q^{\ast}}})}\bigr)\bigr|,
\end{align*}
where we used the stationarity of $(X_i)_{i\in\No}$ in the last line. We have
\begin{align*}
	&J_n(p,q) 
	\\&\leq 
	\sum_{i_q=0}^{n-1} \sum_{\substack{0\leq i_1,\ldots,i_{q-1},i_{q+1},\ldots,i_p\leq i_q \\ i_p^{\ast} \leq n-1}} 
	A_{i_1,\ldots,i_p} 
		+ \sum_{i_q=0}^{n-1} \sum_{\substack{0\leq i_1,\ldots,i_{q-1},i_{q+1},\ldots,i_p\leq i_q \\ i_p^{\ast} \leq n-1}} 
	B_{i_1,\ldots,i_p}.
\end{align*}
An application of the multiple mixing property \eqref{con:multmix_pol} yields
\begin{align*}
	\sum_{i_q=0}^{n-1} \sum_{\substack{0\leq i_1,\ldots,i_{q-1},i_{q+1},\ldots,i_p\leq i_q \\ i_p^{\ast} \leq n-1}} 	A_{i_1,\ldots,i_p} 
	&\leq K\|f(X_0)\|_r \|f\|_{\mathcal{G}} \sum_{i_q=0}^{n-1}   
	(i_q+1)^{p-1} \Theta(i_q) \\
	&\leq K' \|f(X_0)\|_r \|f\|_{\mathcal{G}}
\end{align*}
for some constant $K'<\infty$, since $\sum_{i_q=0}^{\infty} i_q^{p-1}\Theta(i_q)<\infty$ by \eqref{eq:temp022}.
Finally 
\begin{align*}
	&\sum_{i_q=0}^{n-1} \sum_{\substack{0\leq i_1,\ldots,i_{q-1},i_{q+1},\ldots,i_p\leq i_q \\ i_p^{\ast} \leq n-1}} 
	B_{i_1,\ldots,i_p} \\
	&\leq\sum_{i_q=0}^{n-1} \Biggl\{ \sum_{\substack{0\leq i_1,\ldots,i_{q-1}\leq n-1 \\ i_{q-1}^{\ast} \leq n-1}} 
	\Bigl| \E\bigl({f(X_0)f(X_{i_1^{\ast}}) \cdot\ldots\cdot f(X_{i_q^{\ast}})}\bigr) \bigr| \\
	&\hspace{3ex}\cdot
	\sum_{\substack{0\leq i_{q+1},\ldots,i_p\leq n-1 \\ i_p^{\ast}-i_q^{\ast} \leq n-1}}
	\bigl| \E\bigl({f(X_{0})f(X_{i_{q+1}}) \cdot\ldots\cdot f(X_{i_p^{\ast}-i_q^{\ast}})}\bigr) \bigr| \Biggr\} \\
	&= n I_n(p-1)I_n(p-q).
\end{align*}
\qed

\begin{lemma}\label{lem:temp002}
If $(X_n)_{n\in\No}\in\MMP{\Theta}{r}{\mathcal{G}}$, then for all $p\in\N$ such that \eqref{eq:temp022} is satisfied there is a constant $K_p<\infty$, such that 
\begin{align}
 I_n(p) & \leq K_p \sum_{i=1}^{\ugauss{p/2}} n^{i-1} \| f(X_0)\|_r^i \|f\|_{\mathcal{G}}^i \label{eq:temp023}
\end{align}
for all $f\in\mathcal{G}$ with $\|f\|_{\infty}\leq 1$ and $\E(f(X_0))=0$.
\end{lemma}

\proof
We will use mathematical induction to prove the lemma. By Lemma \ref{lem:tmp_bound_Jpq} we can easily see that 
\[ I_n(1) \leq K_1 \| f(X_0) \|_r \: \|f\|_{\mathcal{G}} \]
for some constant $K_1<\infty$ if \eqref{eq:temp022} is satisfied. Now consider an arbitrary $\tilde{p}\geq 2$ satisfying \eqref{eq:temp022} and assume that \eqref{eq:temp023} holds for all $p\leq \tilde{p}-1$. We have
\begin{align*} 
I_n(\tilde{p})
 &\leq \sum_{q=1}^{\tilde{p}} J_n(\tilde{p},q) \\
 &\leq  \sum_{q=1}^{\tilde{p}} \Bigl( K'\| f(X_0)\|_r \|f\|_{\mathcal{G}} + n I_n(q-1)I_n({\tilde{p}}-q) \Bigr)\\
\CUT{
&\leq {\tilde{p}}K' \| f(X_0)\|_r \|f\|_{\mathcal{G}} \\
 &\hspace{3ex} + n \!\sum_{q=1}^{\tilde{p}}\! \Bigl(K_{q-1} \sum_{i=1}^{\ugauss{\frac{q-1}{2}}} n^{i-1} \| f(X_0)\|_r^i \|f\|_{\mathcal{G}}^i\Bigr) \Bigl( K_{\tilde{p}-q} \!\sum_{j=1}^{\ugauss{\frac{\tilde{p}-q}{2}}}\! n^{j-1} \| f(X_0)\|_r^j \|f\|_{\mathcal{G}}^j \Bigr) \\
 }
&\leq K'{\tilde{p}}\| f(X_0)\|_r \|f\|_{\mathcal{G}} + n \sum_{q=1}^{\tilde{p}} K''\sum_{i=2}^{\ugauss{\frac{q-1}{2}}+\ugauss{\frac{\tilde{p}-q}{2}}} n^{i-2} \| f(X_0)\|_r^i \|f\|_{\mathcal{G}}^i \\
\CUT{ &\leq K'' \biggl\{ \| f(X_0)\|_r \|f\|_{\mathcal{G}} + n {\tilde{p}} \sum_{i=2}^{\ugauss{\tilde{p}/2}} n^{i-2} \| f(X_0)\|_r^i \|f\|_{\mathcal{G}}^i\biggr\} \\}
 &\leq K_{\tilde{p}} \sum_{i=1}^{\ugauss{\tilde{p}/2}}n^{i-1} \|  f(X_0)\|_r^i \|f\|_{\mathcal{G}}^i
\end{align*}
for some constants $K',K'',K_{\tilde{p}}<\infty$, since $\ugauss{\frac{q-1}{2}}+\ugauss{\frac{\tilde{p}-q}{2}} \leq \ugauss{\frac{\tilde{p}}{2}}$.
\qed

\proof[Proposition \ref{pro:multimix->2p-moment}]
By \eqref{eq:temp024} and Lemma \ref{lem:temp002} we immediately obtain
\begin{align*}
	\E\Bigl( \Bigl| \sum_{i=1}^n f(X_i) \Bigr|^{2p} \Bigr) \leq (2p)!n I_n(2p-1)   
	\leq K_p \sum_{i=1}^{p} n^{i} \| f(X_0)\|_r^i \|f\|_{\mathcal{G}}^i
\end{align*}
since \eqref{eq:temp020} implies that \eqref{eq:temp022} holds with $p$ replaced by $2p-1$.
\qed


\subsection{Proof of Proposition \ref{pro:causal->multimix}}

Since $(X_i)_{i\in\No}$ is a causal function of an i.i.d.\ process, we have a representation
$$X_i= G(\xi_{i},\xi_{i-1},\ldots),$$
with $G:\mathcal{X}^{\No}\ra\R^d$.
 Let $(\xi'_j)_{j\in\Z}$ and $(\xi''_j)_{j\in\Z}$ be copies of the underlying process $(\xi_j)_{j\in\Z}$ such that all three processes are independent. Set
\begin{align*}
	\dot{X}_i^{(k)} & := G(\xi_{i},\xi_{i-1},\ldots,\xi_{i-k+1},\xi'_{i-k},\xi'_{i-k-1},\ldots) \\
	\ddot{X}_i^{(k)} & := G(\xi_{i},\xi_{i-1},\ldots,\xi_{i-k+1},\xi''_{i-k},\xi''_{i-k-1},\ldots) 
\end{align*}
and note that therefore $(X_i)_{i\in\No} \dequal (\dot{X}_i^{(k)})_{i\in\No} \dequal (\ddot{X}_i^{(k)})_{i\in\No}$. We have
\begin{align}
	&\bigl| \Cov\bigl(f(X_0) \ldots f(X_{i_{q-1}^{*}})\,,\,f(X_{i_{q}^{*}})\ldots f(X_{i_{p}^{*}}) \bigr) \bigr| \notag\\
	&\leq 
	\bigl| \Cov\bigl(f(X_0) \ldots f(X_{i_{q-1}^{*}})-f(\dot{X}^{(k)}_0)\ldots f(\dot{X}^{(k)}_{i_{q-1}^{*}})\,,%
	\,f(X_{i_{q}^{*}})\ldots f(X_{i_{p}^{*}}) \bigr) \bigr| \notag\\
	&\hspace{3ex}+\bigl| \Cov\bigl(f(\dot{X}^{(k)}_0) \ldots f(\dot{X}^{(k)}_{i_{q-1}^{*}})\,,%
	\,f(X_{i_{q}^{*}})\ldots f(X_{i_{p}^{*}}) - f(\ddot{X}^{(k)}_{i_{q}^{*}})\ldots f(\ddot{X}^{(k)}_{i_{p}^{*}}) \bigr) \bigr| \notag\\
	&\hspace{3ex}+\bigl| \Cov\bigl(f(\dot{X}^{(k)}_0) \ldots f(\dot{X}^{(k)}_{i_{q-1}^{*}})\,,%
	\,f(\ddot{X}^{(k)}_{i_{q}^{*}})\ldots f(\ddot{X}^{(k)}_{i_{p}^{*}}) \bigr) \bigr|. \label{eq:temp028}
\end{align}
Since $f(\dot{X}^{(k)}_0) \cdot\ldots\cdot f(\dot{X}^{(k)}_{i_{q-1}^{*}})$ is $\sigma\bigl(\{\xi_j:  j\leq i_{q-1}^{*}\}\cup\{\xi'_{j}: j\in\Z\}\bigr)$-measurable while 
$f(\ddot{X}^{(k)}_{i_{q}^{*}})\cdot\ldots\cdot f(\ddot{X}^{(k)}_{i_{p}^{*}})$ is $\sigma\bigl( \{ \xi_j: j > i_q^*-k \} \cup \{ \xi''_{j}: j\in\Z\}\bigr)$-measurable,
the functions in the last covariance on the right-hand side of \eqref{eq:temp028} are independent as soon as $k\leq i_q$  and thus the last summand is equal to $0$ in this case.

Recall that we only consider such $f$ that satisfy $\|f\|_{\infty}\leq 1$. If we apply H\"older's inequality to equation \eqref{eq:temp028} we obtain for $r,s$ satisfying $\frac{1}{r}+\frac{1}{s}=1$
\begin{align}
	&\bigl|\Cov\bigl(f(X_0) \cdot\ldots\cdot f(X_{i_{q-1}^{*}})-f(\dot{X}^{(k)}_0)\cdot\ldots\cdot f(\dot{X}^{(k)}_{i_{q-1}^{*}})\,,%
	\,f(X_{i_{q}^{*}})\cdot\ldots\cdot f(X_{i_{p}^{*}}) \bigr) \bigr| \notag\\
	&\leq 2 \| f(X_0) \cdot\ldots\cdot f(X_{i_{q-1}^{*}})-f(\dot{X}^{(k)}_0)\cdot\ldots\cdot f(\dot{X}^{(k)}_{i_{q-1}^{*}}) \|_s
	\| f(X_{i_{q}^{*}})\cdot\ldots\cdot f(X_{i_{p}^{*}}) \|_r \notag\\
	&\leq  2q \|f(X_0)\|_r  \| f(X_0)-f(\dot{X}^{(k)}_0)\|_s \label{eq:temp029}
\end{align}
where we used that for $a_i, b_i\in[-1,1]$,
\[\Bigl| \prod_{i=1}^n a_i - \prod_{i=1}^n b_i \Bigl| \leq \sum_{i=1}^n \bigl| a_i-b_i\bigr|. \]
Since $|f(x)-f(y)| \leq \|f\|_{\mathcal{H}_{\alpha}} \|x-y\|^{\alpha}$,
an application of Jensen's inequality to \eqref{eq:temp029} yields
\begin{align*}
	&\bigl|\Cov\bigl(f(X_0) \cdot\ldots\cdot f(X_{i_{q-1}^{*}})-f(\dot{X}^{(k)}_0)\cdot\ldots\cdot f(\dot{X}^{(k)}_{i_{q-1}^{*}})\,,%
	\,f(X_{i_{q}^{*}})\cdot\ldots\cdot f(X_{i_{p}^{*}}) \bigr) \bigr| \notag\\
	&\leq 2q \|f(X_0)\|_r \|f\|_{\mathcal{H}_{\alpha}} \bigl(\| X_0-\dot{X}^{(k)}_0\|_s\bigr)^{\alpha} \notag\\
	&= 2q \|f(X_0)\|_r \|f\|_{\mathcal{H}_{\alpha}}  (\delta_{k,s})^{\alpha}.  
\end{align*}	
Analogously we can show that 
\begin{align*}
	&\bigl| \Cov\bigl(f(\dot{X}^{(k)}_0) \ldots f(\dot{X}^{(k)}_{i_{q-1}^{*}})\,,
	\,f(X_{i_{q}^{*}})\ldots f(X_{i_{p}^{*}}) - f(\ddot{X}^{(k)}_{i_{q}^{*}})\ldots f(\ddot{X}^{(k)}_{i_{p}^{*}}) \bigr) \bigr| \notag\\
	&\leq 2 (p-q) \|f(X_0)\|_r \|f\|_{\mathcal{H}_{\alpha}} \| (\delta_{k,s})^{\alpha},
\end{align*}
thus for $k=i_q$ we have
\begin{align*}
	\bigl| \Cov\bigl(f(X_0) \ldots f(X_{i_{q-1}^{*}})\,,\,f(X_{i_{q}^{*}})\ldots f(X_{i_{p}^{*}}) \bigr) \bigr| 
	&\leq 2q \|f(X_0)\|_r \|f\|_{\mathcal{H}_{\alpha}} (\delta_{i_q,s})^{\alpha}. \notag
\end{align*}
This proves the first part of the proposition. The second part is a consequence of Proposition \ref{pro:multimix->2p-moment}.
\qed

\subsection{Proof of Proposition \ref{pro:clt_causal}}
We will apply the following result obtained by Dedecker in \cite{Ded98}.
\begin{proposition}\label{pro:crit_dedeker}
Let $(Y_i)_{i\in\No}$ be an ergodic stationary process with $\E(Y_0)=0$ and $\E(Y_0^2)<\infty$, which is adapted to a filtration $(\mathcal{M}_i)_{i\in\No}$. If
\[
\sum_{i=0}^n Y_0 \E(Y_i | \mathcal{M}_0)\;\mbox{ converges in }\L_1,
\]
then
\[ \frac{1}{\sqrt{n}} \sum_{i=1}^n Y_i \dconv N(0,\sigma^2)\ \ \text{ as }n\gu,\]
where $\sigma^2=\E(Y_0^2)+2\sum_{i=1}^{\infty} \E(Y_0 \cdot Y_i)<+\infty.$
\end{proposition}

\proof[Proposition \ref{pro:clt_causal}]
Choose an arbitrary $f\in\mathcal{H}_{\alpha}$ with $\E(f(X_0))=0$. The process $(Y_i)_{i\in\No}$ given by
$Y_i:= f(X_i)$
is centered, ergodic, has finite second moments and is adapted to the filtration
\[ (\mathcal{M}_i)_{i\in\No}:=\bigl(\sigma(\xi_i,\xi_{i-1},\ldots)\bigr)_{i\in\No}.\]
As before, let $(\xi'_j)_{j\in\Z}$ be an independent copy of $(\xi_j)_{j\in\Z}$ and set
\begin{align*}
	X'_i&:= G(\xi'_{i},\xi'_{i-1},\ldots), \\
	\dot{X}'_i&:= G(\xi'_{i},\xi'_{i-1},\ldots,\xi'_{1},\xi_{0},\xi_{-1},\ldots).
\end{align*}
Observe that by the independence of $\mathcal{M}_0$ and $\sigma(\{\xi'_i: i\in\Z\})$ we have that
\begin{align*}
	\E\bigl(f(X'_i)|\mathcal{M}_0\bigr) &= \E(f(X'_i)) = 0, &
	\E\bigl(f(X_i)|\mathcal{M}_0\bigr) &= \E\bigl(f(\dot{X}'_i)|\mathcal{M}_0\bigr).
\end{align*}
Thus
\begin{align*}
	\E\Bigl\{\bigl| Y_0 \E\bigl(Y_i|\mathcal{M}_0\bigr) \bigr|\Bigr\} 
	&\leq \|f\|_{\infty} \E\Bigl\{\bigl| \E\bigl(f(X_i)|\mathcal{M}_0\bigr) \bigr|\Bigr\} \allowdisplaybreaks[1]\\
	&= \|f\|_{\infty} \E\Bigl\{\bigl| \E\bigl(f(\dot{X}'_i)-f(X'_i)|\mathcal{M}_0\bigr) \bigr|\Bigr\} \allowdisplaybreaks[1]\\
	&\leq \|f\|_{\infty} \E\bigl| f(\dot{X}'_i) - f(X'_i) \bigr| \allowdisplaybreaks[1]\\
	&\leq \|f\|_{\mathcal{H}_{\alpha}}^2 \E \bigl| \dot{X}'_i - X'_i \bigr|^{\alpha} \allowdisplaybreaks[1]\\
	&\leq \|f\|_{\mathcal{H}_{\alpha}}^2 (\delta_{i,1})^{\alpha},
\end{align*}
where we used Jensen's inequality and $f\in{\mathcal{H}_{\alpha}}$ in the last steps. 
Therefore, by \eqref{con:causal}, $$\sum_{i=1}^n Y_0 \E\bigl(Y_i|\mathcal{M}_0\bigr)$$ converges in $\L_1$ and thus Proposition \ref{pro:crit_dedeker} applies.
\qed

\paragraph*{Acknowledgements}
This article will be published in ``Journal of Theoretical Probability'' (Springer, Heidelberg).
The final publication will be available at springerlink.com.

We would like to thank the anonymous referee for valuable comments and suggestions which helped to improve the first version of the paper.
We are also grateful to Herold Dehling for many helpful discussions.

This research was partially supported by German Research Foundation grant DE 370-4 project: \textit{New Techniques for Empirical Processes of Dependent Data}.


\end{document}